\documentclass[11pt]{article}

\usepackage{setspace}

\usepackage[margin=1in]{geometry}

\usepackage[english]{babel}
\usepackage{tikz}

\usepackage{subfigure}
\usepackage[comma,authoryear]{natbib}  
\if@usehyper
\RequirePackage[colorlinks=false,allbordercolors={1 1 1}]{hyperref}
\fi

\usepackage{mdframed}
\usepackage{amsmath}
\usepackage{graphicx}
\usepackage{amssymb}
\usepackage{mathtools}
\usepackage{amsthm}
\usepackage[colorlinks=true, allcolors=blue]{hyperref}
\usepackage{url} 
\usepackage{subfiles} 
\usepackage{amsfonts}
\usepackage{amsthm}
\usepackage{subcaption}
\usepackage[section]{placeins}
\usepackage{bbm} 
\usepackage{enumitem}
\setlist{nosep}  
\usepackage{lmodern}
\usepackage{textcomp} 

\usepackage{algorithm}
\usepackage{algpseudocode}
\usepackage{bm}  

\usepackage{mathtools}
\DeclarePairedDelimiterX{\norm}[1]{\lVert}{\rVert}{#1}

\newtheorem{theorem}{Theorem}[section]
\newtheorem{lemma}[theorem]{Lemma}
\newtheorem{definition}[theorem]{Definition}
\newtheorem{remark}[theorem]{Remark}

\newcommand{\bd}{\boldsymbol}
\usepackage{caption}

\DeclareMathOperator*{\argmin}{arg\,min}

\usepackage{xcolor}

\definecolor{green}{rgb}{0.0, 0.5, 0.0}

\definecolor{darkgreen}{RGB}{0,100,0}

\usepackage{placeins}

\usepackage[normalem]{ulem}

\usepackage{dsfont}  

\newcounter{assumption}[section]
\renewcommand{\theassumption}{A.\arabic{assumption}} 

\newcommand{\assumption}[1]{\refstepcounter{assumption}\item[(\theassumption)]\label{#1}}

\usepackage[normalem]{ulem} 
\usepackage{hyperref}

\usepackage{authblk}
\title{Linear Quasi-Shrinkage Estimator for High-Dimensional Optimization with Linear Constraints}
\date{}
\author[]{Naqi Huang}
\author[]{Nestor Parolya}
\author[]{Theresia van Essen}
\affil[]{Delft Institute of Applied Mathematics, Delft University of Technology\\

Mekelweg 4, 2628 CD Delft, The Netherlands}

\newcommand{\bA}{\mathbf{A}_m}
\newcommand{\bx}{\mathbf{x}}
\newcommand{\bU}{\mathbf{U}_m}
\newcommand{\bI}{\mathbf{I}_m}
\newcommand{\bbA}{\bar{\mathbf{A}}_m}
\newcommand{\tbA}{\tilde{\mathbf{A}}^k_m}
\newcommand{\bS}{\boldsymbol\Sigma_m}
\newcommand{\bE}{\mathbf{E}_m}
\newcommand{\dotprod}[2]{\left\langle#1,\, #2\right\rangle}

\begin{document}
\maketitle

\begin{abstract}
In large-scale data-driven optimization problems, parameters are often only known approximately due to noisy and small-sized samples. We consider optimization problems with linear constraints where the true parameter matrix is not precisely known, and the number of constraints and variables are both large. Our goal is to construct a linear estimator of the true parameter matrix by minimizing the Frobenius distance between the estimator and the true parameter matrix. Our method offers three key advantages: 1) the coefficients of the linear estimator are consistently estimated from the observations and require no further calibration; 2) it only requires the sample size to be greater than one and it delivers stable performance
across varied sample sizes; and 3) the constraints of the formulated optimization problem using
the estimator remain linear, ensuring computational efficiency when the number of constraints
and variables are large. Simulation shows that our linear estimator consistently produces stable outcomes in terms of the objective value, the ratio of violated constraints and the magnitude of constraint violation across various scenarios, compared to the nominal and robust methods. Additionally, it demonstrates resilience against high levels of noise, making it a robust choice under uncertainty. 
\end{abstract}


\section{Introduction}
Optimization problems with linear constraints have demonstrated their effectiveness as a practical framework for numerous decision problems, given that many real-world problems can be represented or closely approximated using linear constraints. In practice, due to various types of errors such as measurement errors, estimation errors and implementation errors, only noisy samples of the parameters are accessible, which is referred to as `uncertainty'. Many optimization problems involve uncertain parameters, including high-dimensional problems in which the numbers of constraints and variables are both large, with the ratio of the number of constraints and the number of variables bounded by a finite constant. This paper considers optimization problems with linear constraints where the parameters are known only through limited noisy samples, and the numbers of constraints and variables increase to infinity while their ratio converges to a finite constant. 
In our setting, only a small amount of data per parameter is available, a phenomenon frequently observed in real life (see, e.g., \cite{xu2016statistical}). The setting characterized by a small sample size and high-dimensional parameter space is referred to as the small-data, large-scale regime in \cite{gupta2021small}. As described in that context, the combination of features of highly detailed decision making, time-changing environments and low-precision estimates drives the small-data, large-scale phenomenon. Problems such as new-user product recommendations and disaster response operations exhibit these features (see, \cite{gupta2021small}).

Optimization under uncertainty has developed rapidly in recent years. The main approaches to address those optimization problems fall under the scope of robust optimization, stochastic optimization 
 and distributionally robust optimization. 
 The general idea of robust optimization is to define an uncertainty set for the uncertain parameters, and then, optimize against the worst case within this set. While robust optimization protects against worst-case scenarios, it may result in overly conservative solutions since the low violation of constraints comes at the cost of the conservatism of the model (e.g., \cite{ben2000robust}, \cite{bertsimas2004price}, \cite{nguyen2020distributionally} and \cite{bertsimas2022}). In contexts where probabilistic information is available, stochastic or distributionally robust approaches may yield more realistic solutions. Standard stochastic optimization aims to optimize expected performance under a known probability distribution of uncertain parameters (see, e.g., \cite{dantzig1955linear} and \cite{shapiro2021lectures}). A well-known method within the framework of stochastic optimization is chance-constrained optimization, which focuses on satisfying the constraints with a specific probability level (see, e.g., \cite{charnes1959chance} and \cite{shapiro2021lectures}). Distributionally robust optimization combines the ideas of stochastic programming and traditional robust optimization, aiming to optimize worst-case expected performance over an ambiguity set of probability distributions (see, e.g., \cite{delage2010distributionally} and \cite{rahimian2019distributionally}). Over the last decade, the availability of data has increased significantly. There have been data-driven methodologies developed in the aforementioned fields. Sample average approximation (SAA), a well-known data-driven approach in stochastic programming, presents a good performance in the large sample regime (see, e.g., \cite{kleywegt2002sample} and \cite{shapiro2021lectures}). However, small samples may result in poor approximations. An increasing number of works also has focused on data-driven robust optimization (see, e.g., \cite{bertsimas2018data}), including methods on distributional robustness (see, e.g., \cite{delage2010distributionally}, \cite{mohajerin2018data} and \cite{rahimian2022frameworks}). These approaches are generally effective in finite-sample regimes and have demonstrated good performance. Nevertheless, many of them introduce exogenous parameters, such as radii of ambiguity sets that require careful calibration. When the sample sizes are very small, techniques like cross-validation or bootstrapping used for such calibration may become unreliable or impractical (see, e.g., \cite{isaksson2008cross}). Moreover, in the works that adopt hypothesis testing in data-driven robust optimization, the results also depend on the choice of the testing method (see, e.g., \cite{bertsimas2018data}).

Motivated by \cite{xu2016statistical}, we adopt a statistical approach to high-dimensional optimization that avoids the need to calibrate or tune exogenous parameters. In the context of high-dimensional optimization problems with linear constraints, our aim is to estimate the high-dimensional parameter matrix given a small sample size. One can expect that the solution obtained through the estimator of the parameter matrix is close to the true solution if the estimator matrix is a `nice' estimator of the true parameter matrix in the asymptotic sense. In random matrix theory, there has been extensive research on this topic. For example, \cite{lw2004} establishes a well-conditioned estimator for a large dimensional covariance matrix through shrinkage. A non-linear shrinkage method is developed in \cite{ledoit2012nonlinear} and a quadratic shrinkage method in \cite{ledoit2022quadratic}. While estimating a covariance matrix is not required, the concept of the shrinkage method can still be utilized and it has been widely exploited in portfolio selection problems (see, e.g., \cite{golosnoy2007multivariate},  \cite{frahm2010dominating} and \cite{bodnar2022optimal}). Shrinkage has also been studied in optimization beyond portfolio selection. \cite{davarnia2017estimation} propose shrinking a maximum likelihood estimator towards a prespecified vector before using it in an optimization problem. \cite{gupta2022data} develop Shrunken-SAA, which pools data across stochastic optimization problems and learns the amount of pooling in a small-data regime. Our focus is on estimating a noisy rectangular constraint matrix from a fixed number of matrix observations, using consistently estimable scalar quasi-shrinkage coefficients while preserving the linearity of the constraints.


We assume an additive noise model for our parameter matrix, resulting in large information-plus-noise type matrices, more generally classified as large non-centered random matrices. This type of matrices has attracted considerable attention in various applications, for example, in wireless communication (see, e.g.,  \cite{moustakas2003mimo}, \cite{dumont2010capacity}, \cite{hachem2012large} and \cite{hachem2013subspace}). Numerous properties regarding this type of random matrices have also been explored (see, e.g., \cite{banna2020clt} and \cite{zhou2023limiting}). \cite{xu2016statistical} also apply this type of matrices in various optimization problems. 

This paper builds upon the work of \cite{xu2016statistical} who adopt a robust model under an ellipsoidal uncertainty set to address linear optimization problems featuring a high-dimensional parameter matrix. In a similar setting, particularly when having a small sample size, we propose a linear estimator that combines the sample average with a prespecified target matrix. The target matrix can be selected to reflect available structural or prior information about the parameter matrix and the shrinkage coefficients are consistently estimated from the observations, thereby avoiding the tuning and calibration of exogenous parameters required by \cite{xu2016statistical} and other data-driven methods (see, e.g., \cite{mohajerin2018data}, \cite{bertsimas2018data}, \cite{nguyen2020distributionally} and \cite{long2023robust}). This makes our approach particularly suitable for regimes with limited data, where techniques like cross-validation and bootstrapping may become unreliable or difficult to apply. Moreover, this study draws on ideas from \cite{lw2004} and \cite{bodnar2014strong}, extending the linear shrinkage method from focusing on covariance matrices to a broader class of matrices. The term `quasi-shrinkage' adopted in this paper reflects the fact that unlike in classical linear shrinkage,  the coefficients in our setting are not necessarily positive nor should they sum up to one. This arises from the more general nature of the matrices considered,  which offers greater flexibility in the resulting weights.

The contributions of our work are as follows: 
\begin{enumerate}
    \item We propose a linear quasi-shrinkage estimator for the unknown parameter matrix, specifically designed for high-dimensional settings where the numbers of constraints and variables are both large, with the ratio of the number of constraints and the number of variables bounded by a finite constant.
    \item Our method requires a small sample size (only more than one observation of the parameter matrix is needed), compared to other existing data-driven robust approaches.
    \item The target matrix can reflect available structural or prior information about the parameter matrix, while the shrinkage coefficients are consistently estimated from the observations without additional tuning and calibration.
    \item The linearity of the constraints is preserved, making the method computationally efficient, particularly in high-dimensional settings.
    \item We show that the asymptotic loss of our estimator is always lower than the naive `plug-in' approach under stated assumptions.
    \item Simulation results show that our method achieves low constraint violation while maintaining the objective value close to the true optimum.
\end{enumerate}

The structure of this paper is as follows. 
First, we elaborate on the considered problem setting in Section \ref{Sec:Problem Setup}. In Section \ref{Sec:Linear estimation without decision variables}, we propose a linear quasi-shrinkage estimator. More specifically, we assume that the observations of the parameter matrix follow an additive correlated noise model in Section \ref{Sec:Model Formulation} and the linear estimation procedure is then presented in Section \ref{Sec:Linear estimation procedure},
followed by an asymptotic analysis of the estimator in Section \ref{Sec:Analysis} and a more specific discussion on incorporating prior information in Section \ref{Sec:choice of U}. Simulation results and conclusions are presented in Section \ref{Sec:Simulation} and Section \ref{Sec:Conclusion}, respectively.


\section{Problem setup}\label{Sec:Problem Setup}
We address the following general optimization problem:
\begin{equation}\label{true model}
    \max\limits_{\mathbf{x} \in \mathbf{X}} \left\{f(\mathbf{x}): \mathbf{A}\mathbf{x} \leq \mathbf{b}\right\}
\end{equation}
where $\mathbf{A} \in M_{m\times p}(\mathbb{R})$, $\mathbf{X} \subseteq \mathbb{R}^p$ and $\mathbf{b} \in \mathbb{R}^m$. The function $f(\mathbf{x})$ is the objective function to be optimized over $\mathbf{x}$, for which we impose no restriction. We assume that the parameter matrix $\mathbf{A}$ is known solely through noisy samples, hence represents data-driven constraints, and there are $n$ observations for $\mathbf{A}$, denoted by $\tilde{\mathbf{A}}^k$ for $k=1, 2, \dots, n$ ($n > 1$). We focus particularly on the high-dimensional regime, where the numbers of constraints $m$ and variables $p$ are both large, with the ratio of the number of constraints and the number of variables bounded by a finite constant. We specifically focus on the scenario where $m$ and $p$ tend towards infinity such that the ratio $m/p$ approaches a constant $c\in[0, \infty)$, while keeping $n$ small and finite.

Estimating the true matrix $\mathbf{A}$ by simply taking the average of the observations is not reliable: on the one hand, as the sample size $n$ is small, the applicability of the Law of Large Numbers (LLN) is questionable; on the other hand, even if we have a large sample size $n$, the distance between the matrices could still be large if both the number of constraints $m$ and the number of variables $p$ increase. Consequently, as $m$ and $p$ go to infinity, relying on the sample mean of observations for each parameter is not a viable solution. This limitation can be explained by the noise accumulation in the entries of the estimator of the matrix $\mathbf{A}$.

Therefore, we propose a linear quasi-shrinkage method to derive an estimator for the true parameter matrix (see, e.g., \cite{lw2004} and \cite{bodnar2014strong}). The term `quasi' appears  in our framework because the coefficients of the linear combination are not restricted to be positive, nor should they sum up to one. We compare this approach with both the nominal method (plug-in approach) and the data-driven robust method proposed in \cite{mohajerin2018data}. 
The effectiveness of our linear quasi-shrinkage estimator is demonstrated through theoretical findings and simulations.

\section{Linear quasi-shrinkage estimator}\label{Sec:Linear estimation without decision variables}
In this section, we show how the proposed linear quasi-shrinkage estimator of the parameter matrix is derived under certain assumptions. Given that the ratio of the number of constraints $m$ and the number of variables $p$ is bounded by a finite constant, they are implicitly assumed to be dependent on each other. Thus, we introduce the subscript $m$ on relevant matrices to indicate that they are associated with the number of constraints $m$, and hence implicitly with the corresponding number of variables $p$.

\subsection{Assumptions and interpretations}\label{Sec:Model Formulation}
We begin by introducing a set of technical assumptions for our linear quasi-shrinkage method.
\vspace{0.8em}
\begin{enumerate}[label=(\textbf{\theassumption}), ref=\theassumption]
    \assumption{ass:bounded trace} $\sup\limits_{m} \left(\frac1{mp}\mathrm{tr}\left(\mathbf{A}_m\mathbf{A}_m^\top\right)\right) < \infty$.
    \vspace{0.8em}
    \assumption{ass:high dimensional regime} High-dimensional asymptotic regime: $m\rightarrow \infty$, $p=p(m) \rightarrow \infty$ such that  $m/p \rightarrow c \in [0, \infty)$.
    
    \vspace{0.8em}
    \assumption{ass:correlated noise model} There are $n > 1$ observations ($n$ small and finite) for $\mathbf{A}_m$, which can be represented by $\tilde{\mathbf{A}}_m^k = \mathbf{A}_m +  \mathbf{\Sigma}_m^{1/2}\mathbf{E}_m^k$ for $k=1,2,\dots,n$, where $\mathbf{\Sigma}_m$ is an $m\times m$ unknown covariance matrix with the property that $0 < \inf\limits_{m}\left(\frac1m\mathrm{tr}\left(\mathbf{\Sigma}_m\right)\right) \leq \sup\limits_{m}\left(\frac1m\mathrm{tr}\left(\mathbf{\Sigma}_m\right)\right) < \infty$ and $\mathbf{\Sigma}^{1/2}_m$ denotes the symmetric square root matrix of $\mathbf{\Sigma}_m$ that is positive definite. The matrices $\mathbf{E}_m^k$ for $k=1,2,\dots, n$ are mutually independent, and each has i.i.d. entries of mean zero, variance equal to one and uniformly bounded $4+\varepsilon$ moments for some small $\varepsilon > 0$.
 \end{enumerate}
\vspace{0.8em}
Assumption \eqref{ass:bounded trace} is a technical requirement needed in the proof of Theorem \ref{thm:cor-asymptotics}. The high-dimensional regime described in Assumption \eqref{ass:high dimensional regime} is the classical setting in random matrix theory. The matrix $\mathbf{\Sigma}_m^{1/2}\mathbf{E}_m^k$ in Assumption \eqref{ass:correlated noise model} is the noise matrix with covariance matrix $\mathbf{\Sigma}_m$ unknown. This implies that the entries within the same column of the true parameter matrix $\mathbf{A}_m$ are assumed to be dependently perturbed. The finite $4 +\varepsilon$ moments and the bounded normalized traces of $\bS$ in Assumption \eqref{ass:correlated noise model} are required for establishing the asymptotic properties shown in Theorem \ref{thm:cor-asymptotics}. Note that if \,$\bS=\sigma^2\bI$ where $\bI$ is an $m\times m$ identity matrix, the main results can be established with finite second moment using the classical Strong Law of Large Numbers (SLLN). Because of the application of this model in wireless communication and signal processing, it is also called an information-plus-noise type model (see, e.g., \cite{zhou2023limiting}).

The next remark shows that we can handle the row-correlated case in a similar way for $m, p = p(m) \rightarrow \infty$ such that $p/m$ converges to a finite constant.
\begin{remark}
 In the same fashion, assume that every row of $\tilde{\mathbf{A}}^k_m$ is correlated in the same manner by 
 \begin{equation}
     \tilde{\mathbf{A}}_m^k = \mathbf{A}_m + \mathbf{E}_m^k \mathbf{\Sigma}_p^{1/2},
 \end{equation}
 where $\mathbf{\Sigma}_p^{1/2}$ is a $p\times p$ symmetric square root of the row-covariance matrix $\mathbf{\Sigma}_p$. This problem can be transformed to the column-correlated scenario by transposing $\tilde{\mathbf{A}}^k_m$,
 \begin{equation}
  \left(\tilde{\mathbf{A}}_m^k\right)^\top = \mathbf{A}_m^\top + \mathbf{\Sigma}_p^{1/2} \left(\mathbf{E}^k_m\right)^\top,   
 \end{equation}
and the covariance matrix is now $\mathbf{\Sigma}_p$. Therefore, to incorporate this scenario, it would be sufficient to transpose the observation matrices and switch the roles of $m$ and $p$ for $m, p = p(m) \rightarrow \infty$ such that $p/m$ converges to a finite constant.
\end{remark}


\begin{remark}
In line with the findings from \cite{zhou2024analysis}, the noise matrix $\mathbf{E}_m^k$ can be adapted to encompass more generic dependency conditions on its matrix entries. Implementing this in our context is possible, but it would significantly complicate the proofs.
\end{remark}

Note that the sample average obeys the following stochastic representation:
\begin{equation}\label{average of A}
\bar{\mathbf{A}}_m =\frac1n\sum\limits_{k=1}^n \tilde{\mathbf{A}}^k_m = \mathbf{A}_m + \frac{1}{\sqrt{n}}\mathbf{\Sigma}_m^{1/2} \left(\frac1{\sqrt{n}}\sum\limits_{k=1}^n \mathbf{E}_m^k\right) = \mathbf{A}_m +  \frac{1}{\sqrt{n}}\mathbf{\Sigma}_m^{1/2} \mathbf{E}_m\,
\end{equation}
where $\mathbf{E}_m:= \frac1{\sqrt{n}}\sum\limits_{k=1}^n \mathbf{E}_m^k$ and $\mathbf{E}_m$ is again a matrix with i.i.d. entries of mean zero, variance equal to one and finite $4 + \varepsilon$  moments, for some small $\varepsilon > 0$.
As mentioned earlier, our objective is to develop an estimator for the true parameter matrix using observed data. In this context, we directly apply the linear quasi-shrinkage method and consequently, the estimator for the true matrix is formulated as follows:
\begin{equation}\label{A_m^*}
    \mathbf{A}^*_m (\alpha_m,\beta_m) = \alpha_m \bar{\mathbf{A}}_m + \beta_m \mathbf{U}_m,  
\end{equation}
where $\mathbf{U}_m$ is an $m \times p$ matrix specified by the user in advance, referred to as the target matrix. The target matrix is selected to reflect available structural or prior information about the parameter matrix. The idea behind the formulation is to `shrink' the matrix $\bar{\mathbf{A}}_m$ to the target matrix $\mathbf{U}_m$. Naturally, the optimal coefficients $\alpha_m$ and $\beta_m$ for this linear quasi-shrinkage method can be determined by minimizing the distance between true $\mathbf{A}_m$ and $\mathbf{A}^*_m$ based on the chosen loss function. Expressed in mathematical terms, the optimal $\alpha_m$ and $\beta_m$ can be found as follows:
\begin{equation}
    (\alpha^*_m, \beta^*_m) =\argmin\limits_{(\alpha_m, \beta_m)}\left\{d(\mathbf{A}_m, \mathbf{A}_m^*(\alpha_m,\beta_m))\right\}.
\end{equation}
Here, $d(\cdot ,\cdot )$ represents the distance between $\mathbf{A}_m$ and $\mathbf{A}_m^*$. Subsequently, $\alpha^*_m$ and $\beta_m^*$ are dependent on the true matrix $\mathbf{A}_m$, which is unknown. This issue is addressed by estimating $\alpha^*_m$ and $\beta_m^*$ from given data in the high-dimensional regime. The procedure to estimate $\alpha^*_m$ and $\beta_m^*$ is elaborated on in Section \ref{Sec:Linear estimation procedure}. Thus, after obtaining the consistent estimators (see Definition \ref{def: consistent estimator} in Appendix \ref{App: Miscellaneous}) $\hat{\alpha}^*$ and $\hat{\beta}^*$ for $\alpha_m^*$ and $\beta_m^*$, respectively, we plug them into \eqref{A_m^*} and solve the following optimization problem:
\begin{equation}\label{linsh model}
\max\limits_{\mathbf{x} \in \mathbf{X}} \left\{f(\mathbf{x}): \mathbf{A}_m^*(\hat{\alpha}^*, \hat{\beta}^*) \mathbf{x} \leq \mathbf{b}\right\}.
\end{equation}

The formulation can be interpreted from the following perspectives.

\vspace{1em}

\noindent\textbf{Connection to nominal approach.}
The nominal method is a purely plug-in approach in which the parameter matrix is replaced by the sample average, i.e.,
\begin{equation}\label{eq: nominal method}
    \max\limits_{\mathbf{x}\in \mathbf{X}} \left\{f\left(\mathbf{x}\right): \bbA \bx \leq \mathbf{b}\right\}.
\end{equation}
Note that $\bar{\mathbf{A}}_m = \mathbf{A}_m^*  = \alpha_m \bar{\mathbf{A}}_m +\beta_m \mathbf{U}_m$ with $\alpha_m =1$ and $\beta_m =0$.

\vspace{1em}

\noindent\textbf{Connection to robust optimization.}
Based on the concept of uncertainty sets in the field of robust optimization (see \cite{bertsimas2022}), the formulation of \eqref{linsh model} can be interpreted as a robust model for nonnegative coefficients $\alpha_m$ and $\beta_m$ under the following  two additional assumptions:
\vspace{0.8em}
\begin{enumerate}[label=(\textbf{\theassumption}), ref=\theassumption]
    \assumption{assump: non-negativity}  All entries of $m \times p$ matrix $\bar{\mathbf{A}}_m$ and decision variables $\mathbf{x}$ are nonnegative.
    \vspace{0.8em}
    \assumption{assump: all ones}  All entries of the $m \times p$ target matrix $\mathbf{U}_m$ are equal to one.
\end{enumerate}
\vspace{0.8em}
Considering $\tilde{\mathbf{a}}_i^k$ (the $i$-th row of $\tilde{\mathbf{A}}_m^k$) as a noisy copy of the true $\mathbf{a}_i$ (the $i$-th row of $\mathbf{A}_m$), we compute the average of $\tilde{\mathbf{a}}^k_i ~(k = 1,2, \dots, n)$ over the $n$ observations and denote it by $\bar{\mathbf{a}}_i$. Precisely, $\bar{\mathbf{a}}_i$ is the vector constructed from the $i$-th row of $\bar{\mathbf{A}}_m$. In the spirit of robustness, we want our constraints to hold for all values in the uncertainty set specified below, leading to the following formulation:
\begin{subequations}\label{robust interpretation}
    \begin{align}
      & \max  \quad f(\mathbf{x}) \\
        &\ \text{s.t.}\quad \left(y_i\bar{\mathbf{a}}_i + \mathbf{z}_i\right)^\top\mathbf{x} \leq b_i, \quad |y_i|\leq \alpha_m, ||\mathbf{z}_i||_{\infty} \leq \beta_m, \quad i = 1,2, \dots, m, \label{robust interpretation: constraint} \\
        & \quad\quad\quad \mathbf{x}\in \mathbf{X},
    \end{align}
\end{subequations}
where parameters $y_i$ and $\mathbf{z}_i$ are considered as the uncertainty parameters lying in box uncertainty sets (see, e.g., \cite{bertsimas2022}), and $\alpha_m$ and $\beta_m$ are two nonnegative constants associated with robustness. The norm $||\cdot||_{\infty}$ is the infinity norm which is defined as the largest absolute value of the entries of the vector. Note that larger $\alpha_m$ and $\beta_m$ result in a more robust and conservative model. Thus, under Assumptions \eqref{assump: non-negativity} and \eqref{assump: all ones}, we have
\begin{equation}\label{robust interpretation counterpart}
\begin{aligned}
    \max_{|y_i|\leq \alpha_m, \|\mathbf{z}_i\|_{\infty} \leq \beta_m}\ \left(y_i\bar{\mathbf{a}}_i + \mathbf{z}_i\right)^\top\mathbf{x} &=\max_{|y_i| \leq \alpha_m} \ y_i \bar{\mathbf{a}}_i^\top \mathbf{x} + \max_{||\mathbf{z}_i||_{\infty} \leq \beta_m} \mathbf{z}_i^\top\mathbf{x}\\
    &= (\alpha_m \bar{\mathbf{a}}_i)^\top\mathbf{x} + \left(\beta_m\mathbf{1}_{p}\right)^\top \mathbf{x} \\
    &=  (\alpha_m \bar{\mathbf{a}}_i + \beta_m\mathbf{1}_{p})^\top \mathbf{x},
\end{aligned}
\end{equation}
 where ${\bf 1}_{p}$ denotes the $p$-dimensional vector of all ones. This implies that (\ref{robust interpretation: constraint}) is equivalent to $(\alpha_m \bar{\mathbf{a}}_i  +\beta_m{\bf 1}_{p} )^\top\mathbf{x}\leq b_i, i = 1,2, \dots, m$, which is referred to as robust counterpart in the literature (see, e.g., \cite{bertsimas2022}). Therefore, the formulation of \eqref{robust interpretation} is equivalent to 
 \begin{equation}
    \max\limits_{ \mathbf{x} \in \mathbf{X}} \left\{f(\mathbf{x}): (\alpha_m \bbA + \beta_m\bU ) \mathbf{x} \leq \mathbf{b}\right\}
 \end{equation}
 under Assumptions \eqref{assump: non-negativity} and \eqref{assump: all ones}. 

\begin{remark}
    In practical applications, a coefficient preceding an uncertain parameter, e.g., $y_i$ in formulation \eqref{robust interpretation: constraint}, typically represents multiplicative noise (bias). For instance, \cite{ben2000robust} examine multiplicative perturbations, where uncertain parameters are scaled relative to observed values. Additionally, the term $\mathbf{z}_i$ in our formulation captures additive perturbations, which are standard in classical robust optimization frameworks (see, e.g., \cite{bertsimas2022}). Moreover, in control systems and related fields, it is common to model both additive and multiplicative uncertainties jointly (see, e.g., \cite{munoz2015robust}). 
\end{remark}
\noindent\textbf{Bayesian interpretation.}
The Bayesian interpretation of our formulation aligns with the insights that are already discussed in \cite{lw2004}. In our formulation, $\mathbf{U}_m$  represents the prior information we know beforehand and $\bar{\mathbf{A}}_m$ represents the sample information. With the sample information revealed, we shrink sample information $\bar{\mathbf{A}}_m$ toward the target matrix $\mathbf{U}_m$ to integrate our prior beliefs. When no prior information is available, we set all entries of $\mathbf{U}_m$ equal to one, implying that each parameter is treated equally. If we possess complete knowledge of the true parameter matrix $\mathbf{A}_m$, then $\mathbf{U}_m$ can be set to $\mathbf{A}_m$ and one can expect that the optimal $\alpha_m$ and $\beta_m$ are $0$ and $1$, respectively. More details on the target matrix are given in Section \ref{Sec:choice of U}.

\subsection{Linear estimating procedure}\label{Sec:Linear estimation procedure}
As mentioned in the previous section, we aim to minimize the distance between the true parameter matrix $\mathbf{A}_m$ and the estimator $\mathbf{A}_m^*$. Following \cite{leung1987estimation} and \cite{lw2004}, we choose to incorporate the Frobenius norm for the distance between matrices. The former works, however, concentrate their attention on the expected Frobenius distance, while we consider this distance without taking the average, which in fact implies a slightly different loss function. The Frobenius norm is a widely used matrix norm in matrix analysis, which is unitarily invariant and resembles the Euclidean norm for vectors. It is a matrix norm of an $m \times p$ matrix $\mathbf{A}$ defined as the square root of the sum of the squares of its entries (see, e.g., \cite{horn2012matrix}),
which is equal to
\begin{equation}
    ||\mathbf{A}|| = \sqrt{\text{\rm tr}\left(\mathbf{A} \mathbf{A}^\top\right)}.
\end{equation}
For the remainder of the paper, we adopt the following inner product notation: 
\begin{equation}
    \left< \mathbf{A},\,\mathbf{B}\right> :=\mathrm{tr}\left(\mathbf{A}\mathbf{B}^\top \right). 
\end{equation}
Note that $\|\mathbf{A}\|^2 =  \mathrm{tr}\left(\mathbf{A}\mathbf{A}^\top \right) = \left< \mathbf{A},\,\mathbf{A}\right>$. Our goal is to determine the optimal linear combination $\mathbf{A}^*_m (\alpha_m, \beta_m ) = \alpha_m \bar{\mathbf{A}}_m +\beta_m \mathbf{U}_m$ over $\alpha_m$ and $\beta_m$, where the squared Frobenius distance $||\mathbf{A}^*_m- \mathbf{A}_m
||^2$ is minimized. As we are minimizing the Frobenius distance, we take its square for the sake of simplicity. The study by \cite{haff1980empirical} explores this class of linear shrinkage estimators for covariance matrices, but they do not provide any asymptotic optimality results. On the contrary, both \cite{lw2004} and \cite{bodnar2014strong} achieved optimality in the asymptotic sense. Subsequently, they developed `bona fide' estimators for the shrinkage coefficients $\alpha_m$ and $\beta_m$. The concept of `bona fide' estimators emphasizes that, for a fixed target matrix, these estimators are entirely data-driven, independent of any unknown parameters, and do not require calibration procedures. Our approach draws inspiration from \cite{lw2004} and \cite{bodnar2014strong} with the key difference being that we apply the linear shrinkage method to a general rectangular matrix rather than to a covariance matrix.

We follow the following schematic procedure:
\begin{enumerate}
    \item Find the `finite-dimension oracle' (i.e., oracle for finite dimensions $m$ and $p$) estimators $\alpha_m^* (\mathbf{A}_m, \bar{\mathbf{A}}_m)$ and $\beta_m^*(\mathbf{A}_m, \bar{\mathbf{A}}_m)$. Those estimators are optimal for any fixed $m$ and $p$ but not applicable in practice.
    \item Find the asymptotically equivalent quantities for $\alpha_m^*(\mathbf{A}_m, \bar{\mathbf{A}}_m)$ and $\beta^*_m(\mathbf{A}_m, \bar{\mathbf{A}}_m)$ denoted by $\alpha^*(\mathbf{A}_m)$ and $\beta^*(\mathbf{A}_m)$, respectively. The asymptotic equivalents are not applicable as well but can be consistently estimated. One may call them the `asymptotic oracles'.
    \item Derive the consistent estimators of the asymptotic values $\alpha^* (\mathbf{A}_m)$ and $\beta^*(\mathbf{A}_m)$. This leads to the completely data-driven $\hat{\alpha}^*(\bar{\mathbf{A}}_m)$ and $\hat{\beta}^*(\bar{\mathbf{A}}_m)$, which we call `bona fide' estimators. These can be efficiently used in practice.
\end{enumerate}
The procedure can be summarized by the following diagram:\\

\begin{center}
\begin{tikzpicture}
  \node (astarm) at (1,0){ $\begin{array}{c} \left(\alpha^*_m \left(\mathbf{A}_m, \bar{\mathbf{A}}_m\right), \beta^*_m\left(\mathbf{A}_m, \bar{\mathbf{A}}_m\right)\right) \\ \text{`finite-dimension oracle' estimators} \end{array}$};
  
  \node (astar) at (6,4) {$\begin{array}{c} \left(\alpha^*\left(\mathbf{A}_m\right), \beta^*\left(\mathbf{A}_m\right) \right)\\ \text{`asymptotic oracles'} \end{array}$};
  
  \node (hata) at (11,0)  {$\begin{array}{c}\left(\hat{\alpha}^* \left(\bar{\mathbf{A}}_m\right), \hat{\beta}^* \left(\bar{\mathbf{A}}_m\right)\right)\\ \text{`bona fide' estimators}\end{array}$};

  \draw[->] (astarm) --node[midway, sloped, above] {2.  $m\rightarrow \infty$} (astar);
  \draw[->] (astar) --node[midway, sloped, above] {3. Estimated via $\bar{\mathbf{A}}_m$} (hata);
\end{tikzpicture}
\end{center}

Thus, one may expect that the theoretically optimal but inapplicable $\alpha_m^*$ and the applicable $\hat{\alpha}^*$ are equivalent asymptotically. For simplicity of notation, we suppress the dependence on matrices and write $(\alpha_m^*, \beta_m^*)$, $(\alpha^*, \beta^*)$ and $(\hat{\alpha}^*, \hat{\beta}^*)$ instead of $\left(\alpha^*_m \left(\mathbf{A}_m, \bar{\mathbf{A}}_m\right), \beta^*_m\left(\mathbf{A}_m, \bar{\mathbf{A}}_m\right)\right)$, $\left(\alpha^*\left(\mathbf{A}_m\right), \beta^*\left(\mathbf{A}_m\right) \right)$ and $\left(\hat{\alpha}^* \left(\bar{\mathbf{A}}_m\right), \hat{\beta}^* \left(\bar{\mathbf{A}}_m\right)\right)$, respectively, for the remainder of the paper.

\subsubsection{Finite-dimension oracles}\label{Sec:alpha_m^* and beta_m^*}
In this section, our objective is to determine the optimal $\alpha_m$ and $\beta_m$ that minimize the squared Frobenius distance between the estimator matrix $\mathbf{A}_{m}^*(\alpha_m, \beta_m)$ and the true matrix $\mathbf{A}_{m}$. Note that this step does not lead to applicable estimators since these quantities depend on the unknown matrix $\mathbf{A}_m$. But this step is necessary to derive the estimators of the optimal linear quasi-shrinkage coefficients, as will be discussed in Section \ref{Sec:hat_alpha^* and hat_beta^*}. One may refer to the resulting $\alpha^*_m$ and $\beta^*_m$ as `finite-dimension oracles' (or `finite-dimension optimal' estimators), i.e., they are optimal for all fixed $m$ and $p$ but inapplicable in practice. This result is presented in Theorem \ref{fs_optimal}.
\begin{theorem} \label{fs_optimal}
   Consider the optimization problem
   \begin{subequations}\label{distance optimization problem}
       \begin{align}
         \min\limits_{\alpha_m, \beta_m}&~ \left\|\mathbf{A}_{m}^*(\alpha_m, \beta_m) - \mathbf{A}_{m}\right\|^2\\
         \mathrm{s.t.} &\quad \mathbf{A}_m^*(\alpha_m, \beta_m) = \alpha_m \bar{\mathbf{A}}_m  +\beta_m \mathbf{U}_m.
       \end{align}
   \end{subequations}
The optimization problem \eqref{distance optimization problem} has a unique optimal solution, given by
\begin{subequations}\label{alpha_m* and beta_m*}
\begin{align}
\label{alpha_m*}
\alpha^*_{m}  & = \frac{
\dotprod{\bbA}{\bA} \dotprod{\bU}{\bU}
- \dotprod{\bbA}{\bU} \dotprod{\bA}{\bU}
}{
\dotprod{\bbA}{\bbA} \dotprod{\bU}{\bU}
- \dotprod{\bbA}{\bU}^2
},\\
\label{beta_m*}
 \beta^*_{m}  &= \frac{
\dotprod{\bA}{\bU} -\alpha_m^* \dotprod{\bbA}{\bU}
}
{\dotprod{\bU}{\bU}}
 \end{align}
\end{subequations}
where $\bar{\mathbf{A}}_m$ is defined in (\ref{average of A}) and $\mathbf{U}_m$ is  the $m \times p$ target matrix, assumed to be linearly independent of $\bar{\mathbf{A}}_m$. 
\end{theorem}
\begin{proof}
Denote the loss function
\begin{equation}
L(\alpha_m, \beta_m)=\left\| \mathbf{A}_m^*(\alpha_m, \beta_m) - \mathbf{A}_m \right\|^2 = \dotprod{\mathbf{A}_m^*(\alpha_m, \beta_m) - \mathbf{A}_m}{\mathbf{A}_m^*(\alpha_m, \beta_m) - \mathbf{A}_m},
\end{equation}
which is a quadratic function in $\alpha_m$ and $\beta_m$. By completing the square first in $\beta_m$ and then in $\alpha_m$, we obtain
\begin{equation}\label{loss function}
\begin{aligned}
L(& \alpha_m, \beta_m)
= \dotprod{\bU}{\bU}\left( \beta_m 
- \frac{ \dotprod{\bA}{\bU} - \alpha_m \dotprod{\bbA}{\bU} }{ \dotprod{\bU}{\bU} } \right)^2 \\[1ex]
&+ \frac{ \dotprod{\bbA}{\bbA} \dotprod{\bU}{\bU} - \dotprod{\bbA}{\bU}^2 }{\dotprod{\bU}{\bU} }
\left( \alpha_m 
- \frac{ \dotprod{\bbA}{\bA} \dotprod{\bU}{\bU} - \dotprod{\bbA}{\bU} \dotprod{\bA}{\bU} }{ \dotprod{\bbA}{\bbA} \dotprod{\bU}{\bU} - \dotprod{\bbA}{\bU}^2 } \right)^2 \\[1ex]
&- \left[
\left( 
\frac{ \dotprod{\bbA}{\bA} \dotprod{\bU}{\bU} - \dotprod{\bbA}{\bU} \dotprod{\bA}{\bU} }{ \dotprod{\bbA}{\bbA} \dotprod{\bU}{\bU} - \dotprod{\bbA}{\bU}^2 } 
\right)^2
\frac{ \dotprod{\bbA}{\bbA} \dotprod{\bU}{\bU} - \dotprod{\bbA}{\bU}^2 }{ \dotprod{\bU}{\bU} } \right. \\[1ex]
& \left. \quad + \frac{ \dotprod{\bA}{\bU}^2 }{ \dotprod{\bU}{\bU} } - \dotprod{\bA}{\bA}
\right].
\end{aligned}
\end{equation}
Note that $\dotprod{\bU}{\bU} > 0$ and $\dotprod{\bbA}{\bbA} \dotprod{\bU}{\bU} - \dotprod{\bbA}{\bU}^2 \geq 0$ due to the Cauchy-Schwarz inequality with equality holding if and only if $\bbA$ and $\bU$ are linearly dependent. The loss function $L(\alpha_m, \beta_m)$ is therefore uniquely minimized at the values $\alpha_m^*$ and $\beta_m^*$, given in \eqref{alpha_m*} and \eqref{beta_m*}, respectively. \qedhere
\end{proof}

\begin{remark}
    In practice, $\bU$ is selected in advance by the decision maker, which allows considerable flexibility to ensure that $\bbA$ and $\bU$ are linearly independent. When $\bU$ is intentionally selected to be linearly dependent on $\bbA$, the formulation simplifies to
    \begin{equation}
        \mathbf{A}_m^* = \alpha_m \bbA 
    \end{equation}
where $\mathbf{A}_m^*$ is simply a scalar multiple of $\bbA$. In this case, minimizing $\|\mathbf{A}_m^* -\bA\|^2$ yields $\alpha_m^* = \frac{\dotprod{\bbA}{\bA}}{\dotprod{\bbA}{\bbA}}$. Accordingly, one can derive the corresponding `asymptotic oracle' and `bona fide' estimator for $\alpha_m^*$ as for \eqref{alpha_m* and beta_m*}. Since the derivation follows exactly the same procedure and to maintain consistency and clarity in the main text, we leave this special case out and assume that $\bbA$ and $\bU$ are linearly independent for the remainder of the paper. 


\end{remark}

\begin{remark}\label{tr(bbU) = 1}
    Note that if we consider $\alpha_m^*$ and $\beta_m^*$ as functions of $\mathbf{U}_m$, then
\begin{equation}
    \alpha_m^*(\delta_m \mathbf{U}_m) = \alpha_m^* (\mathbf{U}_m) \quad \beta_m^*(\delta_m \mathbf{U}_m) 
 = \frac1{\delta_m}\beta_m^*(\mathbf{U}_m),
\end{equation}
for any $\delta_m \neq 0$, which leads to 
\begin{equation}
\begin{aligned}
    \mathbf{A}_m^* (\alpha_m^*(\delta_m\mathbf{U}_m), \beta_m^*(\delta_m \mathbf{U}_m)  ) &= \alpha_m^*(\delta_m\mathbf{U}_m) \bar{\mathbf{A}}_m + \beta_m^*(\delta_m \mathbf{U}_m) (\delta_m\mathbf{U}_m)\\
    & = \alpha_m^*(\mathbf{U}_m) \bar{\mathbf{A}}_m + \frac{1}{\delta_m}\beta_m^*(\mathbf{U}_m) (\delta_m\mathbf{U}_m) \\
    & =  \alpha_m^*(\mathbf{U}_m) \bar{\mathbf{A}}_m + \beta_m^*(\mathbf{U}_m) \mathbf{U}_m \\
    & = \mathbf{A}_m^*(\alpha_m^*(\mathbf{U}_m), \beta_m^*(\mathbf{U}_m)).
\end{aligned}
\end{equation}
This means the estimator $\mathbf{A}_m^*$ is independent of the magnitude of $\mathbf{U}_m$. Thus, without loss of generality, we assume that $0 < \inf\limits_{m}\left( \frac1{mp}\dotprod{\bU}{\bU}\right) \leq \sup\limits_{m}\left( \frac1{mp}\dotprod{\bU}{\bU}\right) <\infty$ for the remainder of the paper.


\end{remark}

As already mentioned earlier, the expressions for $\alpha^*_{m}$ and $\beta^*_{m}$ in \eqref{alpha_m*} and \eqref{beta_m*} cannot be directly used in practice in the current form because they contain the unknown true matrix $\mathbf{A}_m$. Therefore, our approach is to derive their asymptotic deterministic equivalents, which are defined as follows.
\begin{definition}[Asymptotic deterministic equivalent]\label{def:Asymptotic deterministic equivalent}
    Suppose we have a sequence of random variables $\tilde{X}_{m,p}$ such that for $m/p \rightarrow c \in [0, \infty)$ as $m \rightarrow \infty, p=p(m)\rightarrow \infty$, we have
    \begin{equation}
        \left|\tilde{X}_{m,p} - X_{m,p}\right| \xrightarrow{a.s.} 0,
    \end{equation}
then the deterministic sequence $X_{m,p}$ is called the asymptotic deterministic equivalent of $\tilde{X}_{m,p}$. We introduce the following notation for this relation:
\begin{equation}
 \tilde{X}_{m,p} \overset{a.e.}{\sim} X_{m,p}.   
\end{equation}

 \end{definition}
 
    Thus, our goal is firstly to find the asymptotic deterministic equivalents for $\alpha_m^*$ and $\beta_m^*$, denoted by $\alpha^*$ and $\beta^*$, which are estimated consistently in the subsequent step as discussed in \cite{bodnar2014strong}. To achieve this, we begin by dividing both the numerator and denominator of $\alpha_m^*$ by $(mp)^2$ in \eqref{alpha_m*}  and those of $\beta_m^*$ by $mp$ in \eqref{beta_m*} to ensure the boundedness of all traces, from which the three terms $\frac{ \dotprod{ \bbA}{ \bU} }{mp},\quad 
\frac{ \dotprod{\bbA}{\bA} }{mp}$, and
$\frac{ \dotprod{\bbA }{\bbA} }{mp}$ naturally arise. Then, we identify the asymptotic deterministic equivalents of those three terms, which directly yield the asymptotic deterministic equivalents of $\alpha^*_m$ and $\beta^*_m$, denoted by $\alpha^*$ and $\beta^*$, respectively. Thereafter, the consistent estimators of $\alpha^*$ and $\beta^*$, denoted by $\hat{\alpha}^*$ and $\hat{\beta}^*$, respectively, can be straightforwardly derived.

\subsubsection{Asymptotic equivalents of \texorpdfstring{$\alpha^*_m $}{alpha} and \texorpdfstring{$\beta^*_{m}$}{beta}}\label{Sec:alpha^* and beta^*}
After the finite-dimension optimal $\alpha^*_m$ and $\beta^*_m$ are found, the next step is to explore their asymptotic behavior as $m, p(m) \to\infty$ such that $m/p\to c \in [0, \infty)$. As stated earlier, it is enough to consider three functionals $\frac{ \dotprod{ \bbA}{ \bU} }{mp},\,
\frac{ \dotprod{\bbA}{\bA} }{mp}$ and
$\frac{ \dotprod{\bbA }{\bbA} }{mp}$. This result is presented in the following theorem.


\begin{theorem}\label{thm:cor-asymptotics}
Under assumptions \eqref{ass:bounded trace}, \eqref{ass:high dimensional regime} and \eqref{ass:correlated noise model}, it holds that
\begin{subequations}\label{traces asym}
\begin{gather}
\frac{1}{mp} \dotprod{\bbA}{\bU} \overset{a.e.}{\sim} \frac{1}{mp} \dotprod{\bA}{\bU}, \label{traces limit: b_AU} \\[1ex]
\frac{1}{mp} \dotprod{\bbA}{\bA} \overset{a.e.}{\sim} \frac{1}{mp} \dotprod{\bA}{\bA}, \label{traces asym: b_AA} \\[1ex]
\frac{1}{mp} \dotprod{\bbA}{\bbA} \overset{a.e.}{\sim} \frac{1}{mp} \dotprod{\bA}{\bA} + \frac{1}{nm} \dotprod{\mathbf{\Sigma}_m}{\bI}. \label{traces asym: b_Ab_A}
\end{gather}
\end{subequations}
\end{theorem}

\begin{proof}
For $\frac{1}{mp} \dotprod{\bbA}{\bU}$, we have
\begin{equation}\label{col b_AU}
\begin{aligned}
\frac{1}{mp} \left| \dotprod{\bbA}{\bU} - \dotprod{\bA}{\bU} \right| 
= \frac{1}{\sqrt{n}mp} \left|\dotprod{ \mathbf{\Sigma}_m^{1/2} \mathbf{E}_m }{ \bU }\right| \xrightarrow{a.s.} 0
\end{aligned}
\end{equation}
and for $\frac{1}{mp} \dotprod{\bbA}{\bA}$, we have
\begin{equation}\label{barA A with Sigma}
\frac{1}{mp} \left| \dotprod{\bbA}{\bA} - \dotprod{\bA}{\bA} \right| 
= \frac{1}{\sqrt{n}mp} \left|\dotprod{ \mathbf{\Sigma}_m^{1/2} \mathbf{E}_m }{ \bA }\right| \xrightarrow{a.s.} 0
\end{equation}
for $m, p(m) \rightarrow \infty$ such that $m/p \rightarrow c \in [0, \infty)$ (see Lemma \ref{APP:correlated KSLLN of A, U with Sigma} in Appendix \ref{app:supplement proofs}).

For $\frac{1}{mp} \dotprod{\bbA}{\bbA}$, we decompose:
\begin{equation}\label{split Abar_Abar with Sigma}
\begin{aligned}
\frac{1}{mp} \dotprod{\bbA}{\bbA}
&= \frac{1}{mp} \dotprod{\bA}{\bA}
+ \frac{2}{\sqrt{n}mp} \dotprod{ \bA }{ \mathbf{\Sigma}_m^{1/2} \mathbf{E}_m } 
+ \frac{1}{nmp} \dotprod{ \mathbf{\Sigma}_m^{1/2} \mathbf{E}_m }{ \mathbf{\Sigma}_m^{1/2} \mathbf{E}_m }\\
& = \frac{1}{mp} \dotprod{\bA}{\bA}
+ \frac{2}{\sqrt{n}mp} \dotprod{ \bA }{ \mathbf{\Sigma}_m^{1/2} \mathbf{E}_m } 
+ \frac{1}{nmp} \dotprod{ \bS}{\bE\bE^\top}
\end{aligned}
\end{equation}
where $\sup\limits_m \left( \frac{1}{mp} \dotprod{\bA}{\bA} \right) < \infty$ by assumption \eqref{ass:bounded trace} and the second term vanishes almost surely as shown in \eqref{barA A with Sigma}. For the third term, $\frac{1}{nmp} \dotprod{ \bS}{\bE\bE^\top}$, we have that
\begin{equation}
    T_m: = \frac{1}{nmp} \dotprod{ \bS}{\bE\bE^\top} = \frac1{np}\sum\limits_{j=1}^p \frac{\mathbf{e}_j^\top \bS \mathbf{e}_j}m
\end{equation}
where $\mathbf{e}_j$ is the $j$-th column vector of $\mathbf{E}_m$. By Assumption \eqref{ass:correlated noise model}, we have that $\mathbb{E}\left(\mathbf{e}_j^\top \bS \mathbf{e}_j\right) =\mathrm{tr}\left(\bS\right)$, and therefore, $\mathbb{E}\left(T_m\right) = \frac1{nm}\mathrm{tr}\left(\bS\right)$. According to Lemma B.26 in \cite{bai2010spectral}, we have
\begin{equation}
   \mathbb{E}\left(\left|\mathbf{e}_j^\top \bS \mathbf{e}_j - \mathrm{tr}\left(\bS\right)\right|^{2+\frac{\varepsilon}2}\right) \leq C\left(\left(\nu_4\mathrm{tr}\left(\bS^2\right)\right)^{1 + \frac{\varepsilon}4} +\nu_{4+\varepsilon} \left(\mathrm{tr}\left(\bS^2\right)\right)^{1+\frac{\varepsilon}4}  \right)  = O\left(\left(\mathrm{tr}\left(\bS\right)\right)^{2 + \frac{\varepsilon}2}\right) 
\end{equation}
where the last equality follows from the fact that $\mathrm{tr}\left(\bS^k\right) \leq \left(\mathrm{tr}\left(\bS\right)\right)^k$ for any $k\geq 1$. Here, $C$ is a constant dependent on $2+\frac\varepsilon2$ only and $\nu_k$ is the $k$-th moment of the entries of $\mathbf{E}_m$. Hence, by the Rosenthal inequality (see, e.g., \cite{rosenthal1970subspaces}), we have that
\begin{equation}
\begin{aligned}
    &\quad\mathbb{E}\left(\left|T_m - \frac{\mathrm{tr}\left(\bS\right)}{nm}\right|^{2+\frac\varepsilon2}\right)  =  \mathbb{E}\left(\left| \frac1{nmp}\sum\limits_{j=1}^p\left(\mathbf{e}_j^\top \bS \mathbf{e}_j -  \mathrm{tr}\left(\bS\right)\right)\right|^{2+\frac\varepsilon2}\right) \\
    &\leq \frac{C}{(nmp)^{2+\frac\varepsilon2}}\left[\left(\sum\limits_{j=1}^p\mathbb{E}\left(\left|\mathbf{e}_j^\top \bS \mathbf{e}_j -  \mathrm{tr}\left(\bS\right) \right|^2\right)\right)^{1+\frac{\varepsilon}4} + \sum\limits_{j=1}^{p}\mathbb{E}\left(\left|\mathbf{e}_j^\top \bS \mathbf{e}_j -  \mathrm{tr}\left(\bS\right)\right|^{2+\frac{\varepsilon}2}\right) \right] \\
    &= \frac{C}{n^{2+\frac\varepsilon2}p^{1+\frac\varepsilon4}} O\left(\left(\frac{\mathrm{tr}\left(\bS\right)}m\right)^{2 +\frac{\varepsilon}2}\right)
\end{aligned}
\end{equation}
where $\sup\limits_{m}\frac{\mathrm{tr}\left(\bS\right)}m < \infty$ under Assumption \eqref{ass:correlated noise model}. Then, by Borel Cantelli lemma (see, e.g. \cite{shiryaev2016probability}), for any $\delta > 0$, Markov's inequality gives
\begin{equation}
\begin{aligned}
    \mathbb{P}\left(\left|T_m - \frac{\mathrm{tr}\left(\bS\right)}{nm}\right| > \frac1{m^{\delta}}\right) &\leq C m^{\delta(2+\frac{\varepsilon}2)} \mathbb{E}\left(\left|T_m - \frac{\mathrm{tr}\left(\bS\right)}{nm}\right|^{2+\frac\varepsilon2}\right)\\ 
    &= O\left(\frac{m^{\delta(2+\frac\varepsilon2)}}{p^{1+\frac{\varepsilon}4}}\right) = O\left(\frac1{m^{1+\frac{\varepsilon}4 - \delta(2+\frac\varepsilon2)}}\right)
\end{aligned}
\end{equation}
where the last equality follows from the fact that $m =O(p)$ under Assumption \eqref{ass:high dimensional regime}. Thus, for $1+\frac{\varepsilon}4 - \delta(2+\frac\varepsilon2) >1$, i.e., $\delta < \frac12 - \frac{2}{4 + \varepsilon}$,
\begin{equation}
    T_m = \frac{1}{nmp} \dotprod{ \bS}{\bE\bE^\top} \overset{a.e.}{\sim} \frac1{nm}\mathrm{tr}\left(\bS\right) = \frac1{nm} \dotprod{\bS}{\bI}.
\end{equation}
Therefore, recalling \eqref{split Abar_Abar with Sigma}, \eqref{traces asym: b_Ab_A} holds. \qedhere

\end{proof}
Replacing the normalized trace terms in $\alpha_m^*$ and $\beta_m^*$ (see \eqref{alpha_m* and beta_m*}) with their asymptotic equivalents derived in Theorem \ref{thm:cor-asymptotics} leads to the following asymptotic deterministic equivalents for $\alpha_m^*$ and $\beta_m^*$:
\begin{subequations}\label{alpha^* and beta^*}
\begin{align}
\label{alpha^*}
\alpha^* &= \frac{\dotprod{\bA}{\bA} \dotprod{\bU}{\bU} - \dotprod{\bA}{\bU}^2}{\left( \dotprod{\bA}{\bA} + \frac{p}{n} \dotprod{\bS}{\bI} \right) \dotprod{\bU}{\bU} - \dotprod{\bA}{\bU}^2}\nonumber\\
& = 1 - \frac{ \frac{p}{n} \dotprod{\bS}{\bI} \dotprod{\bU}{\bU} }
{ \left( \dotprod{\bA}{\bA} + \frac{p}{n} \dotprod{\bS}{\bI} \right) \dotprod{\bU}{\bU} - \dotprod{\bA}{\bU}^2 },\\
\label{beta^*}
\beta^* 
&= \frac{ \dotprod{\bA}{\bU} - \alpha^* \dotprod{\bA}{\bU} }{ \dotprod{\bU}{\bU} }= (1 - \alpha^*) \frac{ \dotprod{\bA}{\bU} }{ \dotprod{\bU}{\bU} }.
\end{align}
\end{subequations}  
Note that the denominator of $\alpha^*$ normalized by $(mp)^2$ gives
\begin{equation}\label{eq: bounded denominator}
\begin{aligned}
    D_m : &= \left(\frac{\dotprod{\bA}{\bA}}{mp} + \frac{\mathrm{tr}\left(\bS\right)}{nm}\right)\frac{\dotprod{\bU}{\bU}}{mp} - \left(\frac{\dotprod{\bA}{\bU}}{mp}\right)^2 \\
    &\geq \frac{\mathrm{tr}\left(\bS\right)}{nm}\frac{\dotprod{\bU}{\bU}}{mp} \geq \frac1n \inf\limits_{m} \frac{\mathrm{tr}\left(\bS\right)}{m} \frac{\dotprod{\bU}{\bU}}{mp}  > 0
\end{aligned}
\end{equation}
where $\inf\limits_{m} \frac{\mathrm{tr}\left(\bS\right)}{m} > 0$ according to Assumption \eqref{ass:correlated noise model} and $\inf\limits_{m}\frac{\dotprod{\bU}{\bU}}{mp} > 0$ according to Remark \ref{tr(bbU) = 1}.
Then, applying Theorem \ref{thm:cor-asymptotics}, we have
\begin{equation}\label{limit alpha^*, beta^*}
\alpha^*_m \overset{a.e.}{\sim} \alpha^*, \quad  \beta^*_m\overset{a.e.}{\sim}\beta^*.
\end{equation}
In \eqref{alpha^* and beta^*}, note that $\dotprod{\bA}{\bA} \dotprod{\bU}{\bU} - \dotprod{\bA}{\bU}^2 \geq 0$ due to Cauchy-Schwarz inequality. Moreover, we have $\dotprod{\bS}{\bI}> 0$ since $\bS$ is symmetric positive definite. It follows that $\alpha^* \in [0,1)$. 
\subsubsection{Consistent estimators of \texorpdfstring{$\alpha^*$}{alpha1} and \texorpdfstring{$\beta^*$}{beta1}}\label{Sec:hat_alpha^* and hat_beta^*}
Up to this point, we have found the asymptotic deterministic equivalents $\alpha^*$ and $\beta^*$ (see \eqref{alpha^* and beta^*}) for $\alpha^*_m$ and $\beta_m^*$ (see \eqref{alpha_m* and beta_m*}). However, they depend on the true parameter matrix $\mathbf{A}_m$ and the covariance matrix $\bd\Sigma_m$ which are unknown. We therefore proceed in this section to implement our next step: estimate $\alpha^*$ and $\beta^*$ consistently in accordance with Definition \ref{def:Asymptotic deterministic equivalent}. This consists of two parts: first, we estimate the two terms: $\frac{1}{mp} \dotprod{\bA}{\bU}$ and $\frac{1}{mp} \dotprod{\bA}{\bA}$ in the formulas of $\alpha^*$ and $\beta^*$; we then estimate the term $\frac1m\dotprod{\bS}{\bI}$. 

The first step is easily accomplished by applying Theorem \ref{thm:cor-asymptotics}: the consistent estimators of \\ $\frac{1}{mp} \dotprod{\bA}{\bU}$ and $\frac{1}{mp} \dotprod{\bA}{\bA}$ 
are given by 
$\frac{1}{mp} \dotprod{\bbA}{\bU}$ and 
$\frac{1}{mp} \dotprod{\bbA}{\bbA} - \frac{1}{nm} \dotprod{\bS}{\bI}$, respectively. To estimate the term $\frac1m\dotprod{\bS}{\bI}$, we introduce the following lemma.

\begin{lemma}\label{lemma: estimator of trace of Sigma}
For any fixed sample size $n>1$, it holds that
\begin{equation}\label{estimator of trace of Sigma}
      \widehat{\sigma}^2_m \overset{a.e.}{\sim} \frac{1}{m} \dotprod{\bS}{\bI}
\end{equation}
where $\widehat{\sigma}^2_m: = \frac{1}{(n-1)mp} \sum\limits_{k=1}^n \dotprod{ \tilde{\mathbf{A}}^k_m - \bbA }{ \tilde{\mathbf{A}}^k_m - \bbA }.$ 
\end{lemma}

\begin{proof}
Note that
\begin{align}
&\frac{1}{(n-1)mp} \sum\limits_{k=1}^n\dotprod{ \tbA - \bbA }{ \tbA - \bbA } \nonumber \\
=~& \frac{1}{n - 1} \sum_{k=1}^n \frac{1}{mp} \dotprod{\tbA}{\tbA}
- \frac{2}{(n - 1)mp} \sum_{k=1}^n \dotprod{\tbA}{\bbA}
+ \frac{n}{(n - 1)mp} \dotprod{\bbA}{\bbA} \nonumber\\
=~& \frac{1}{n - 1} \sum_{k=1}^n \frac{1}{mp} \dotprod{\tbA}{\tbA}
- \frac{2n}{(n - 1)mp} \dotprod{\underbrace{\frac1n\sum\limits_{k=1}^n \tbA}_{=\bbA}}{\bbA}
+ \frac{n}{(n - 1)mp} \dotprod{\bbA}{\bbA}    \nonumber\\
=~&\frac{1}{n - 1} \sum_{k=1}^n \frac{1}{mp} \dotprod{\tbA}{\tbA}
- \frac{n}{n - 1} \cdot \frac{1}{mp} \dotprod{\bbA}{\bbA}. \label{estimator of sigma^2 expansion}
\end{align} 
From \eqref{traces asym: b_Ab_A}, we know that
\begin{equation}\label{traces asym: b_Ab_A_2}
\frac{1}{mp} \dotprod{\bbA}{\bbA}
\overset{a.e.}{\sim}\left( \frac{1}{mp} \dotprod{\bA}{\bA} + \frac{1}{nm} \dotprod{\bS}{\bI}\right). 
\end{equation}
For each $k = 1,2,\dots,n$, comparing $\frac{1}{mp} \dotprod{\tbA}{\tbA}$ with $\frac{1}{mp} \dotprod{\bbA}{\bbA}$, we note that $\bar{\mathbf{A}}_m$ 
equals $\tilde{\mathbf{A}}_m^k$ when $\tilde{\mathbf{A}}^k_m$ is the only observation. Thus,
\begin{equation}\label{traces asym: b_Abarb_Abar}
\frac{1}{mp} \dotprod{\tbA}{\tbA}
\overset{a.e.}{\sim}\left( \frac{1}{mp} \dotprod{\bA}{\bA} + \frac{1}{m} \dotprod{\bS}{\bI} \right).  
\end{equation}
Therefore, when substituting the asymptotic equivalents given in \eqref{traces asym: b_Ab_A_2} and \eqref{traces asym: b_Abarb_Abar} in (\ref{estimator of sigma^2 expansion}), we obtain \eqref{estimator of trace of Sigma}. \qedhere

\end{proof}

From Lemma \ref{lemma: estimator of trace of Sigma}, we conclude that $\widehat{\sigma}^2_m$ is a consistent estimator of $\frac1m\dotprod{\bS}{\bI}$. The proof of Lemma \ref{lemma: estimator of trace of Sigma} also explains why it is necessary to have $n >1$. Since $n$ is a fixed constant, scaling by $\frac1{n}$ instead of $\frac1{n-1}$ would lead to a biased limit: $\frac{n-1}{nm}\dotprod{\bS}{\bI}$. Moreover, when $n=1$, $\tilde{\mathbf{A}}_m^k$ is simply equal to $\bar{\mathbf{A}}_m$, resulting in a degenerate (zero) estimator.

Following \eqref{alpha^* and beta^*} and Lemma \ref{lemma: estimator of trace of Sigma}, we obtain the consistent estimators for $\alpha^*$ and $\beta^*$ given by

\begin{subequations}\label{hat_alpha^* and hat_beta^*} 
\begin{align}
 \hat{\alpha}^* 
&=   1 - \frac{ 
\frac{1}{n(n - 1)} \sum\limits_{k=1}^n \dotprod{\tbA - \bbA}{\tbA - \bbA}  \dotprod{\bU}{\bU} 
}{
\dotprod{\bbA}{\bbA} \dotprod{\bU}{\bU} - \dotprod{\bbA}{\bU}^2
}
\label{hat_alpha^*}, \\
    \hat{\beta}^{*} & = (1 - \hat{\alpha}^*) \frac{ \dotprod{\bbA}{\bU} }{ \dotprod{\bU}{\bU} } \label{hat_beta^*}.
\end{align}
\end{subequations}
The estimators $\hat{\alpha}^*$ and $\hat{\beta}^*$ are now completely data-driven and can be used in practice. They do not require any additional calibration, compared to other data-driven robust optimization methods (e.g., \cite{bertsimas2018data}, \cite{mohajerin2018data}). In the end, we obtain our estimator for the true matrix $\mathbf{A}_m$: 
\begin{equation}\label{bona fide A}
    \hat{\mathbf{A}}_m^* = \hat{\alpha}^*\bar{\mathbf{A}}_m + \hat{\beta}^* \mathbf{U}_m,
\end{equation}
which is then used as input to solve optimization problem \eqref{linsh model}. 

\subsection{Analysis of loss}\label{Sec:Analysis}


In this section, we analyze the asymptotic loss of $\mathbf{A}_m^*=\alpha^*\bbA +\beta^*\bU$. 
Recall from equation \eqref{loss function} that, when substituting the optimal values $\alpha_m^*$ and $\beta_m^*$ into the loss function, i.e., $L^*_m:=L(\alpha^*_m, \beta^*_m) $, the resulting loss for $\bA^*$ is given by
\begin{equation}\label{loss}
\begin{aligned}
L^*_{m} = &- \left[
\left( 
\frac{ \dotprod{\bbA}{\bA} \dotprod{\bU}{\bU} - \dotprod{\bbA}{\bU} \dotprod{\bA}{\bU} }{ \dotprod{\bbA}{\bbA} \dotprod{\bU}{\bU} - \dotprod{\bbA}{\bU}^2 } 
\right)^2
\frac{ \dotprod{\bbA}{\bbA} \dotprod{\bU}{\bU} - \dotprod{\bbA}{\bU}^2 }{ \dotprod{\bU}{\bU} } \right. \\[1ex]
& \left. \quad + \frac{ \dotprod{\bA}{\bU}^2 }{ \dotprod{\bU}{\bU} } - \dotprod{\bA}{\bA}
\right].
\end{aligned}
\end{equation}

\begin{theorem}\label{theorem: asymptotic loss}
Under the conditions of Theorem \ref{thm:cor-asymptotics}, it holds that
\begin{equation}\label{asymptotic loss}
 \frac1{mp}L^*_{m}  \overset{a.e.}{\sim} \frac{\alpha^*}{nm}\dotprod{\bd\Sigma_m}{\mathbf{I}_m}.
 \end{equation}
In other words, $\frac1{mp}\|\bA^*-\bA\|^2  \overset{a.e.}{\sim} \frac{\alpha^*}{nm}\dotprod{\bd\Sigma_m}{\mathbf{I}_m}$ where $\bA^* = \alpha^* \bbA +  \beta^* \bU$. 
\begin{proof}
    See the proof in Appendix \ref{App: proof of asymptotic loss}.
\end{proof}
\end{theorem}

\begin{remark}\label{remark:asymptotic loss from hat A}
    By the consistency of \(\hat\alpha^*\) and \(\hat\beta^*\), together with the boundedness of the normalized matrix norms under Assumption \eqref{ass:bounded trace} and Remark \ref{tr(bbU) = 1}, we have that

\begin{equation}
    \frac{\left\| \hat{\mathbf{A}}_m^* -  \mathbf{A}_m^*\right\|}{\sqrt{mp}} \leq \left|\hat{\alpha}^* - \alpha^*\right|\frac{\|\bbA\|}{\sqrt{mp}} + \left|\hat{\beta}^* - \beta^*\right|\frac{\left\|\bU\right\|}{\sqrt{mp}} \xrightarrow{a.s.} 0
\end{equation}
where $\hat{\mathbf{A}}_m^* = \hat{\alpha}^*\bbA +\hat{\beta}^*\bU$ and $\mathbf{A}_m^* = \alpha^*\bbA +\beta^*\bU$.
Thus, the normalized loss results established in Theorem \ref{theorem: asymptotic loss} also apply to the bona fide estimator \(\hat {\mathbf{A}}_m^*\).
\end{remark}

  Theorem \ref{theorem: asymptotic loss} shows that $\mathbf{A}^*_m$ is not a consistent estimator of $\bA$ in general, as the Frobenius distance between $\mathbf{A}^*_m$ and $\mathbf{A}_m$ does not approach zero asymptotically.

When we consider the loss for the nominal case, i.e., simply plugging in the sample average $\bbA$ as the estimator, we get
\begin{equation}
    \frac{1}{mp}\|\bbA-\bA\|^2 = \frac{1}{mp}\left(\dotprod{\bbA}{\bbA}  -2 \dotprod{\bbA}{\bA} +\dotprod{\bA}{\bA}\right).
\end{equation}
Applying the asymptotic results of Theorem \ref{thm:cor-asymptotics}, we have that
\begin{equation}
     \frac{1}{mp}\|\bbA-\bA\|^2 \overset{a.e.}{\sim}\frac1{nm} \dotprod{\bd\Sigma_m}{\bI}.
\end{equation}
Now, since $\alpha^*$ lies in $[0, 1)$ in our setting, the proposed estimator has a smaller asymptotic normalized Frobenius loss than the nominal estimator. 


\begin{remark}\label{dependence of A and U}
Recall from \eqref{alpha^*} that $\alpha^* =0$ if and only if $\bA$ and $\bU$ are linearly dependent. In this case, with \eqref{beta^*} we obtain
\begin{equation}
    \bA^* = \alpha^* \bbA + \beta^*\bU = \frac{ \dotprod{\bA}{\bU} }{ \dotprod{\bU}{\bU} } \bU = \bA.
\end{equation}
It shows that when $\bA$ and $\bU$ are linearly dependent, the linear quasi-shrinkage estimator exactly recovers $\bA$, resulting in zero asymptotic loss.
\end{remark}
Note that larger $n$ leads to smaller $\frac1{nm}\dotprod{\bd\Sigma_m}{\mathbf{I}_m}$, which in turn increases $\alpha^*$, according to the first equation of \eqref{alpha^*}. This is intuitive, since when $n$ increases, the sample average $\bbA$ becomes more accurate for the true matrix $\bA$, and thus more weight is assigned to $\bbA$. As $n \rightarrow \infty$, $\frac1{nm}$ $\dotprod{\bd\Sigma_m}{\mathbf{I}_m} \rightarrow 0$, which implies $\alpha^* \rightarrow 1$ by \eqref{alpha^*} and consequently, $\beta^* \rightarrow 0$ by \eqref{beta^*}. Hence, the asymptotic loss $\frac1{mp}L_m^*$ approaches zero. This reflects the fact that with a large sample size $n$, the sample average becomes a reliable estimator of the true parameter matrix, as a consequence of the LLN. In contrast, when $n$ is small, $\alpha^*$ is reduced, placing more weight on the target matrix. The resulting weights reflect a trade-off between the quality of the sample average $\bbA$ and the target matrix $\bU$, highlighting the shrinkage behavior of the estimator.

The preceding results concern matrix estimation error. To relate this analysis to the optimization problem, we next
bound the average magnitude of constraint violation in terms
of the Frobenius estimation error and the norm of the resulting
decision vector.
\begin{remark}[Connection to constraint violation]
Let $\hat{\mathbf{x}}$ be a feasible solution for
\begin{equation*}
\max\limits_{\bx \in \mathbf{X}}\left\{\mathbf{c}^\top\bx: \hat{\mathbf{A}}_m^*\bx\leq \mathbf{b}\right\},
\end{equation*}
then $$\mathbf{a}_i^\top\hat{\mathbf{x}}-b_i  = \left(\mathbf{a}_i - \hat{\mathbf{a}}_i^*\right) ^\top\hat{\mathbf{x}}+ \left(\hat{\mathbf{a}}_i^*\right)^\top \hat{\bx}-b_i \leq \left(\mathbf{a}_i - \hat{\mathbf{a}}^*_i\right) ^\top\hat{\mathbf{x}}$$ where $\mathbf{a}_i$ and $\hat{\mathbf{a}}_i^*$ are the $i$-th row of $\mathbf{A}_m$ and $\hat{\mathbf{A}}_m^*$, respectively and $b_i$ is the $i$-th entry of $\mathbf{b}$. Thus, the average magnitude of constraint violation satisfies
\[
\frac{1}{m}\sum_{i=1}^m
\max\left(\mathbf{a}_i^\top\hat{\mathbf{x}}-b_i, 0\right)
\leq
\frac{\|\hat{\mathbf{x}}\|}{\sqrt{m}}
\bigl\|\hat{\mathbf{A}}_m^*-\mathbf{A}_m\bigr\|
\]
where $\frac1{\sqrt{mp}}\bigl\|\hat{\mathbf{A}}_m^*-\mathbf{A}_m\bigr\| \overset{a.e.}{\sim} \sqrt{\frac{\alpha^*}{nm}\dotprod{\bS}{\bI}}$ according to Theorem \ref{theorem: asymptotic loss} and Remark \ref{remark:asymptotic loss from hat A}. 
\end{remark}

\subsection{Choice of target matrix}\label{Sec:choice of U}
As mentioned in Section \ref{Sec:Model Formulation}, $\mathbf{U}_m$ represents prior information of the true parameter matrix. We also briefly discussed in Remark \ref{dependence of A and U} how the asymptotic loss of the estimator can be affected by the choice of $\bU$ in a specific case. In this section, we provide practical guidance on the choice of $\mathbf{U}_m$ given certain types of information of $\mathbf{A}_m$.

\begin{itemize}
    \item \textbf{No prior information on $\mathbf{A}_m$.} When no specific structural information about $\bA$ is available, a simple choice for $\mathbf{U}_m$ is the matrix of all ones, i.e., $\mathbf{U}_m = \mathbf{1}_m\mathbf{1}_p^\top$ where $\mathbf{1}_m$ and $\mathbf{1}_p$ are $m$-dimensional and $p$-dimensional vectors of all ones, respectively, which treats all entries equally. We plug $\hat{\alpha}^*$ and $\hat{\beta}^*$ given by \eqref{hat_alpha^* and hat_beta^*} into \eqref{A_m^*} under this scenario, which gives

\begin{equation*}
\begin{aligned}
\mathbf{A}^*_m & = \hat{\alpha}^* \bar{\mathbf{A}}_m + (1-\hat{\alpha}^*)\frac{\dotprod{\bar{\mathbf{A}}_m }{\mathbf{1}_m\mathbf{1}_p^\top}}{mp} \mathbf{1}_m\mathbf{1}_p^\top=\hat{\alpha}^* \bar{\mathbf{A}}_m + (1-\hat{\alpha}^*)\frac{\mathbf{1}_m^\top\bar{\mathbf{A}}_m \mathbf{1}_p}{mp} \mathbf{1}_m\mathbf{1}_p^\top. 
\end{aligned}
\end{equation*}
Note that $\mathbf{1}_m^\top\bar{\mathbf{A}}_m \mathbf{1}_p$ is the sum of all entries of $\bar{\mathbf{A}}_m$. Therefore, we can see that, without any prior information, every entry of $\bar{\mathbf{A}}_m$ is shrunk towards $\frac{\mathbf{1}_m^\top\bar{\mathbf{A}}_m \mathbf{1}_p}{mp} $, which is simply the average of all entries of $\bar{\mathbf{A}}_m$. This form of the matrix $\bU$ is used further in the simulation study.
    
\item \textbf{Information on scaled $\mathbf{A}_m$ is known}. Suppose the true matrix $\mathbf{A}_m$ is known up to an unknown scaling constant. In this case,  $\bU$ can be taken as a scalar multiple of $\bA$, implying that they are linearly dependent. As discussed in Remark \ref{dependence of A and U}, $\mathbf{A}^*_m$ exactly recovers the true $\bA$ asymptotically under this setting.

\item \textbf{Some entries of $\mathbf{A}_m$ are possibly scaled}.
In the linear quasi-shrinkage method, $\alpha_m$ represents the scale parameter, while $\beta_m$ serves as the shift parameter. Imagine a situation where the decision maker knows that certain elements of the true matrix $\mathbf{A}_m$ may only be scaled but not shifted. In these cases, the corresponding entries in $\mathbf{U}_m$ can be set to zero.

\end{itemize}
Similarly, other scenarios can be considered, for example, when certain entries of the true matrix $\mathbf{A}_m$ are precisely known. In such cases, it is necessary not only to set $\bU$ accordingly but also adapt the noise model, which may alter the entire estimation procedure. We therefore leave this intriguing scenario for future research. Nonetheless, the key takeaway is that incorporating external knowledge about the matrix $\mathbf{A}_m$ into $\mathbf{U}_m$ can potentially improve the performance of the estimator. Furthermore, the linear quasi-shrinkage method can also be extended from a single target matrix to multiple target matrices for greater flexibility (see Appendix \ref{App: multiple targets}).




\section{Simulation}\label{Sec:Simulation}

In this section, we present simulation results demonstrating the effectiveness of our linear quasi-shrinkage method. For simplicity, we consider a linear objective function in our simulations. The optimization problem considered in our simulation is formulated as:
\begin{equation}\label{random problem model}
     \max\limits_{\mathbf{x} \geq \mathbf{0}}\{\mathbf{c}^\top \mathbf{x}:\, \mathbf{A}_m  \mathbf{x}\leq \mathbf{b} \}
\end{equation}
where $\mathbf{c}$ and $\mathbf{b}$ are known but the true matrix $\bA$ is unknown. Instead, we observe a noisy sample of $\bA$, denoted by $\tilde{\mathbf{A}}_m^k$ for $k \in \{1,2,\dots,n\}$, which can be treated as the realization of a random matrix $\tilde{\mathbf{A}}_m$. According to Assumption \eqref{ass:correlated noise model}, we can conclude that true $\mathbf{A}_m$ can be considered as the expectation of $\tilde{\mathbf{A}}_m$. Thus, the true optimization problem then becomes
\begin{equation}
       \max\limits_{\mathbf{x}\geq \mathbf{0}}\left\{\mathbf{c}^\top \mathbf{x}:\, \mathbb{E}_{\mathbb{P}^*}\left(\tilde{\mathbf{A}}_m\right)    \mathbf{x}\leq \mathbf{b}\right\}
\end{equation}
where $\mathbb{P}^*$ is the joint distribution over all entries of $\tilde{\mathbf{A}}_m$. We compare our linear quasi-shrinkage method to the following two methods.

\vspace{0.7em}
\noindent\textbf{Nominal method.} This approach is a naive one, which simply replaces  $\mathbb{E}_{\mathbb{P}^*}\left(\tilde{\mathbf{A}}_m\right)$ by the sample average, i.e., approximating the true distribution by the empirical distribution. The resulting formulation is,
\begin{equation}
    \max\limits_{\mathbf{x}\geq \mathbf{0}}\left\{\mathbf{c}^\top \mathbf{x}:\, \bbA  \mathbf{x}\leq \mathbf{b}\right\}
\end{equation}
where $\bbA: = \frac1n\sum\limits_{k=1}^n \tilde{\textbf{A}}_m^k$  is the empirical average of the observed sample. 

\vspace{0.7em}
\noindent\textbf{Distributionally robust optimization (DRO).}
DRO aims to safeguard against distributional uncertainty by optimizing under the worst-case distribution within a predefined ambiguity set centered at the empirical distribution. Specifically, we formulate the problem using a Wasserstein ball as the ambiguity set (see \cite{mohajerin2018data}):

\begin{equation}\label{dro}
    \max\limits_{\mathbf{x}\geq \mathbf{0}}\left\{\mathbf{c}^\top\mathbf{x}: \sup\limits_{\mathbb{P}_i\in \mathcal{B}_i(\hat{\mathbb{P}}_i, \gamma)}\left(\mathbb{E}_{\mathbb{P}_i}\left(\tilde{\mathbf{a}}_i^\top\mathbf{x}\right)\right) \leq b_i, i =1,2,\dots,m\right\} 
\end{equation}
where $\mathcal{B}_i(\hat{\mathbb{P}}_i, \gamma) := \left\{\mathbb{P}: \|\mathbb{P} -  \hat{\mathbb{P}}_i\|_W \leq \gamma \right\}$ is the Wasserstein ball of radius $\gamma > 0$ around the empirical distribution $\hat{\mathbb{P}}_i$ over $\mathbf{a}_i$, $\mathbf{a}_i$ is the $i$-th row of $\mathbf{A}_m$ and $b_i$ is the i-th entry of $\mathbf{b}$. $\|\cdot\|_W$ is the 1-Wasserstein distance defined in Appendix \ref{App: Miscellaneous},  Theorem \ref{theorem: wasserstein}. Increasing $\gamma$ makes the DRO formulation more conservative,  potentially reducing constraint violations. Note that the worst-case expectation can be reformulated as follows (see proof in Theorem \ref{theorem: wasserstein}):  
\begin{equation}
  \sup\limits_{\mathbb{P}_i\in \mathcal{B}_i(\hat{\mathbb{P}}_i, \gamma)}\left(\mathbb{E}_{\mathbb{P}_i}\left(\tilde{\mathbf{a}}_i^\top\mathbf{x}\right)\right) = \mathbb{E}_{\hat{\mathbb{P}}_i} (\tilde{\mathbf{a}}_i^\top\bx)  + \gamma \|\mathbf{x}\| = \frac1n\sum\limits_{k=1}^n \left(\tilde{\mathbf{a}}_i^k\right)^\top \bx + \gamma \|\mathbf{x}\| = \bar{\mathbf{a}}_i^\top \bx + \gamma \|\mathbf{x}\|,
\end{equation}
where $\tilde{\mathbf{a}}_i^k$  and $\bar{\mathbf{a}}_i$ are the i-th row of $\tilde{\mathbf{A}}_m^k$ and $\bar{\mathbf{A}}_m$, respectively. Therefore, \eqref{dro} can be rewritten as
\begin{equation}\label{dro model}
    \max\limits_{\mathbf{x}\geq \mathbf{0}}\left\{\mathbf{c}^\top\mathbf{x}:  \bar{\mathbf{a}}_i^\top\mathbf{x} + \gamma \|\mathbf{x}\|  \leq b_i,\, i=1,2,\dots,m \right\}.
\end{equation}
Moreover, optimization problem \eqref{dro model} is equivalent to the following robust optimization problem under an ellipsoidal uncertainty set (see, e.g.,
\cite{bertsimas2022}), 
\begin{equation}
    \max\limits_{\mathbf{x}\geq \mathbf{0}}\left\{\mathbf{c}^\top\mathbf{x}:  \left(\bar{\mathbf{a}}_i + \mathbf{z}\right)^\top~\mathbf{x}  \leq b_i,  \|\mathbf{z}\|\leq \gamma, \, i=1,2,\dots, m\right\}.
\end{equation}
In summary, we use DRO model \eqref{dro model} as a benchmark for evaluating our method. The radius $\gamma$ of the Wasserstein ball is treated as an exogenous parameter and is selected using cross-validation (see Appendix \ref{App:cross-validation} for a detailed description).

For the simulations, we draw the true parameter matrix $\mathbf{A}_m$, the right-hand side of constraints $\mathbf{b}$, and the coefficient vector $\mathbf{c}$ in the objective function, all independently from uniform distribution $U(4,6)$.  Assuming that $\bd\Sigma_m=\sigma^2\bd I_m$, each entry of $\mathbf{A}_m$ is then perturbed with i.i.d. Gaussian noise $N(0,\sigma^2)$. We take $\bU$ as the matrix of all ones throughout the simulations, which is well suited to the present setting, where the true matrix entries are independently drawn from the same distribution. Given that the applicable sample size is  $n>1$,  we set the sample size to $n \in \{3, 5\} $ in main comparisons. This experiment is conducted at various noise levels $\sigma \in \{1, 2, 3\}$. Given that this paper focuses on settings where the ratio of the number of constraints $m$ and the number of variables $p$ is bounded by a finite constant, we conduct the simulation with $c = m/p \in \{0.5, 1, 1.5\}$ with $p$ varying from $100$ to $500$ with steps of $100$ in main comparisons. A more detailed investigation with additional values of $c$ and sample size $n$ has also been conducted. For each parameter setting, $50$ repetitions are performed for the three methods: the nominal method, the robust method and the linear quasi-shrinkage method. The performance is assessed across four criteria, as outlined in \cite{xu2016statistical}. These criteria are:
\begin{enumerate}
     \item {\bf Relative objective value:} Defined as the relative error between the objective value obtained by a given method and the objective value obtained under the true constraints, computed as the difference divided by the true objective value.
    \item {\bf Ratio of violated constraints:} Computed by dividing the number of violated constraints by the total number of constraints $m$.
    \item {\bf Magnitude of constraint violation:} Calculated as the average violation across all constraints, computed by summing the amount by which each constraint is violated (when violated) and dividing by the total number of constraints $m$.
    \item {\bf Computation time:} The time required to solve the optimization problem for a single repetition and parameter setting, measured in seconds.
\end{enumerate}

For each of these four criteria, the average result over the 50 repetitions is determined for each of the three methods. Simulation results are shown in the figures in the subsequent pages. All computations were performed using Python 3.13.12 on the DelftBlue high-performance computing cluster at Delft University of Technology. We used Gurobi 12.0.3 to solve the optimization problems. 

\begin{figure}
\begin{subfigure}
\centering
\includegraphics[width=1.05\textwidth]{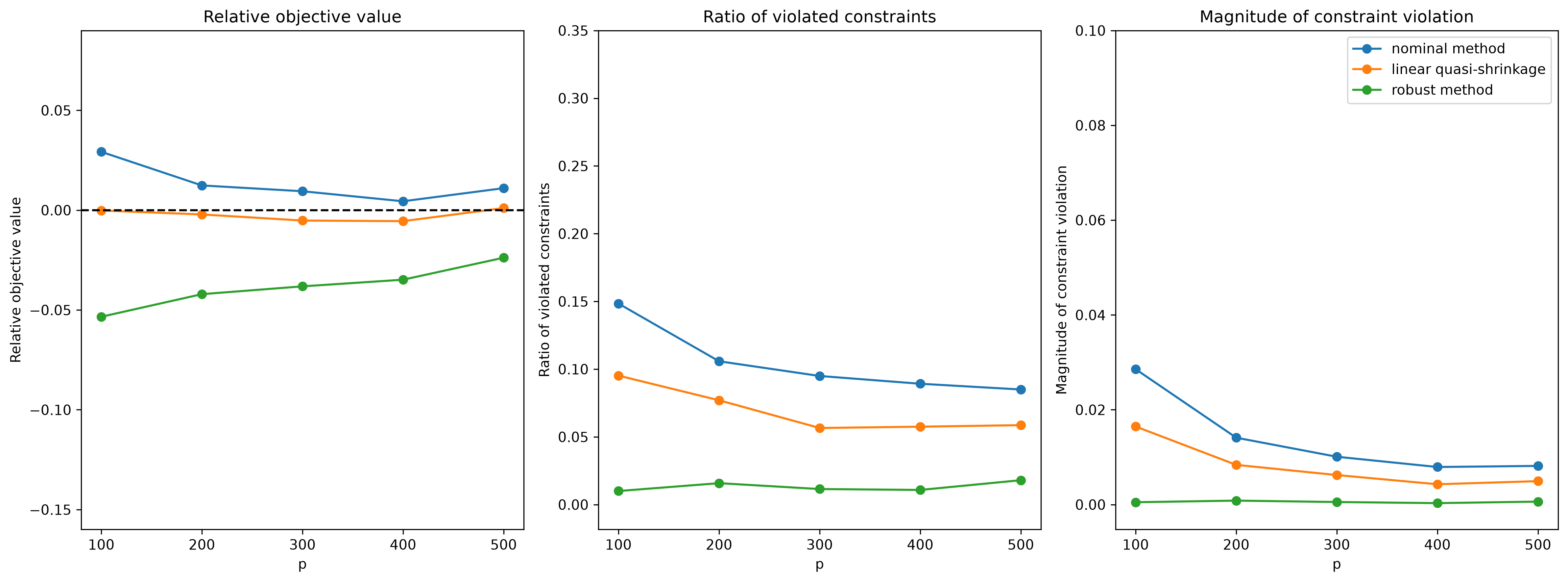}
\caption{Performance for $c=0.5,\, n =3, \, \sigma = 1$}
\label{c05_n3_s1}
\end{subfigure}

\begin{subfigure}
\centering

\includegraphics[width=1.05\textwidth]{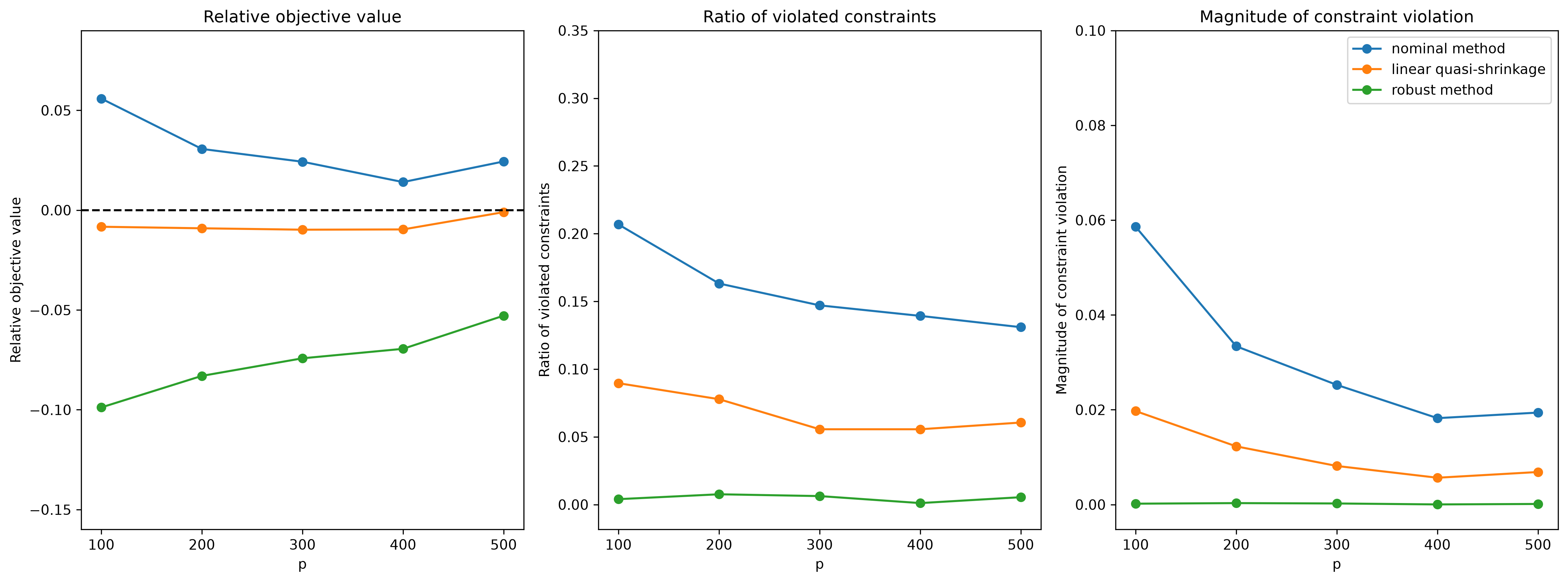}
\caption{Performance for $c=0.5,\, n =3,\, \sigma = 2$}
\label{c05_n3_s2}
\end{subfigure}

\begin{subfigure}
\centering

\includegraphics[width=1.05\textwidth]{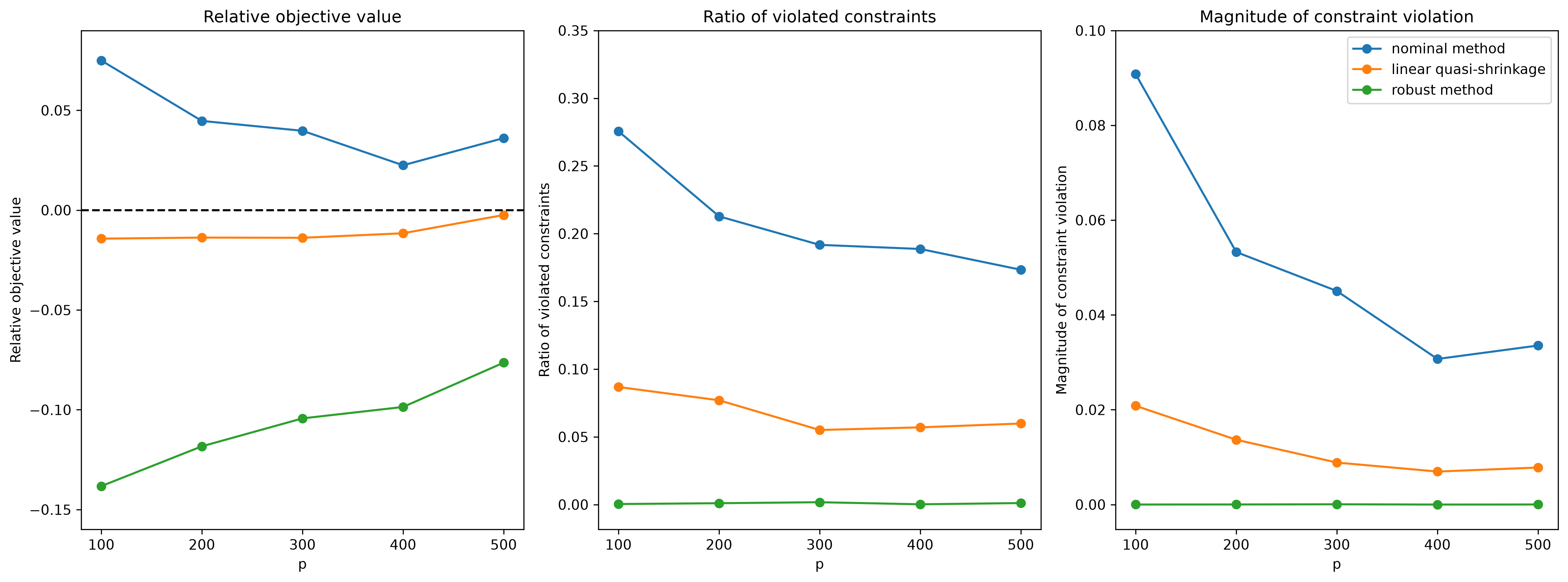}
\caption{Performance for $c=0.5,\, n =3,\, \sigma = 3$}
\label{c05_n3_s3}
\end{subfigure}
\vspace{1.5em}
\end{figure}

In Figures \ref{c05_n3_s1}, \ref{c05_n3_s2} and \ref{c05_n3_s3} with $n=3$, $c=0.5$ and $\sigma = 1, 2,3$, respectively, we compare the performance of three methods: the nominal method (blue), the linear quasi-shrinkage method (orange) and the robust method (green). The comparison is made across three criteria: (1) relative objective value (first column); (2) ratio of violated constraints (second column), and (3) magnitude of constraint violation (third column). In each figure, we vary the number of variables $p$ from $100$ to $500$ with steps of $100$. Based on results, we have the following observations:
\begin{itemize}
    \item In all three figures with varied $\sigma$, increasing $p$ (and accordingly, increasing $m$) generally reduces both the ratio of violated constraints and the magnitude of constraint violation for the nominal and linear quasi-shrinkage methods, although these improvements are not always monotonic. The relative objective values of the linear quasi-shrinkage method remain close to zero, whereas those of the nominal method exhibit larger fluctuations. Meanwhile, the robust method’s relative objective values improve from substantially below zero toward zero with a noticeable remaining gap as $p$ increases.
    
    \item The nominal method generally exhibits a larger ratio of violated constraints and magnitude of constraint violation than the linear quasi-shrinkage method, particularly at higher noise levels.  In contrast, the robust method consistently maintains a low ratio of violated constraints and magnitude of constraint violation, however, at the expense of the objective value. Our linear quasi-shrinkage method strikes a balance between the objective value and constraint violation (i.e., the ratio of violated constraints and the magnitude of constraint violation). In other words, its objective value sticks close to the true objective value while maintaining a relatively low constraint violation.  
        \item A comparison across the three figures shows that as $\sigma$ increases, both the ratio of violated constraints and the magnitude of constraint violation increase noticeably for the nominal method. The linear quasi-shrinkage method maintains a stable ratio of violated constraints and magnitude of constraint violation, while its relative objective values remain close to zero. The robust method exhibits fewer and smaller constraint violations as $\sigma$ increases, but its relative objective values become substantially more negative, indicating a larger relative loss in objective value.
\end{itemize}

\begin{figure}

\begin{subfigure}
\centering

\includegraphics[width=1.05\textwidth]{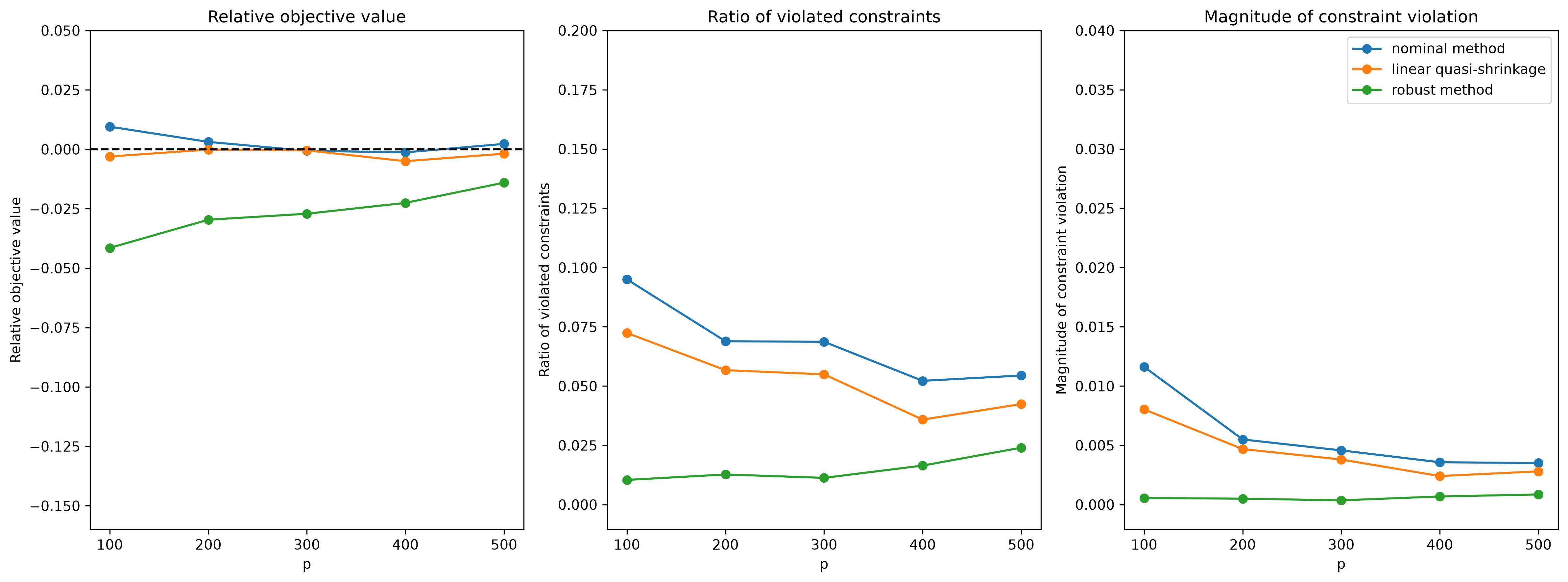}
\caption{Performance for $c=1,\, n =3,\, \sigma = 1$}
\label{c10_n3_s1}
\end{subfigure}
  
\begin{subfigure}
\centering

\includegraphics[width=1.05\textwidth]{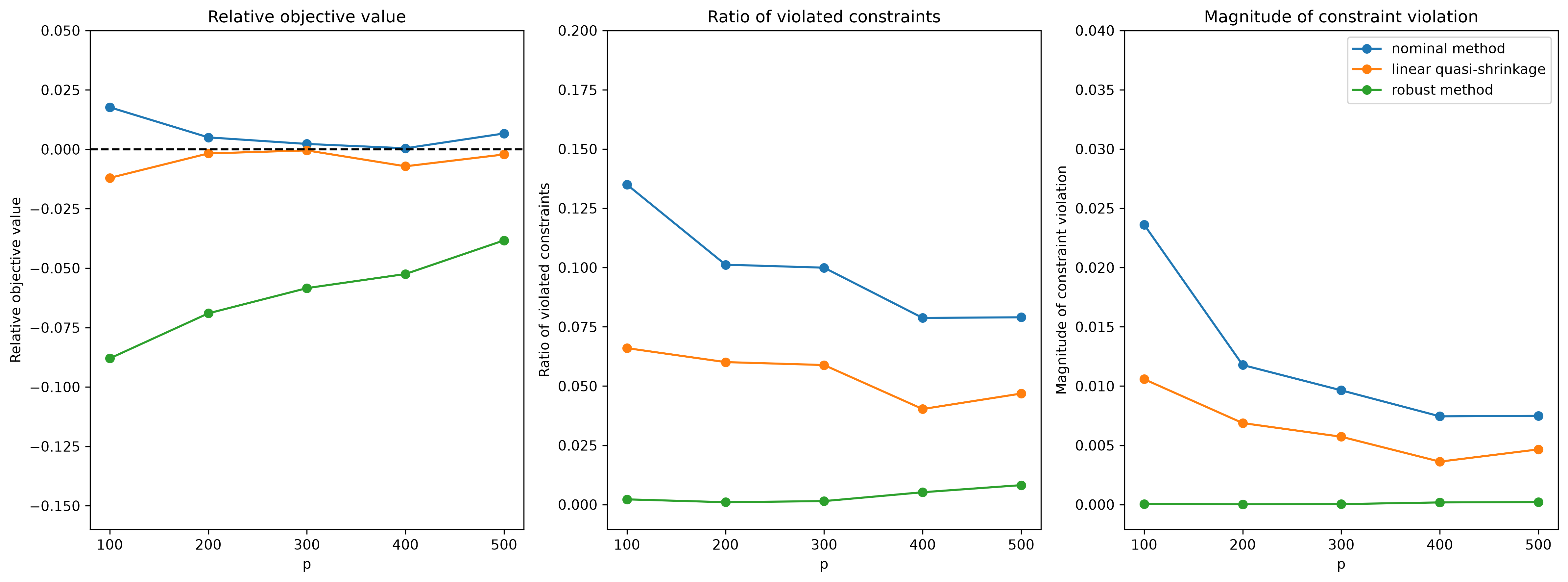}
\caption{Performance for $c=1,\, n =3,\, \sigma = 2$}
\label{c10_n3_s2}
\end{subfigure}
\begin{subfigure}
\centering

\includegraphics[width=1.05\textwidth]{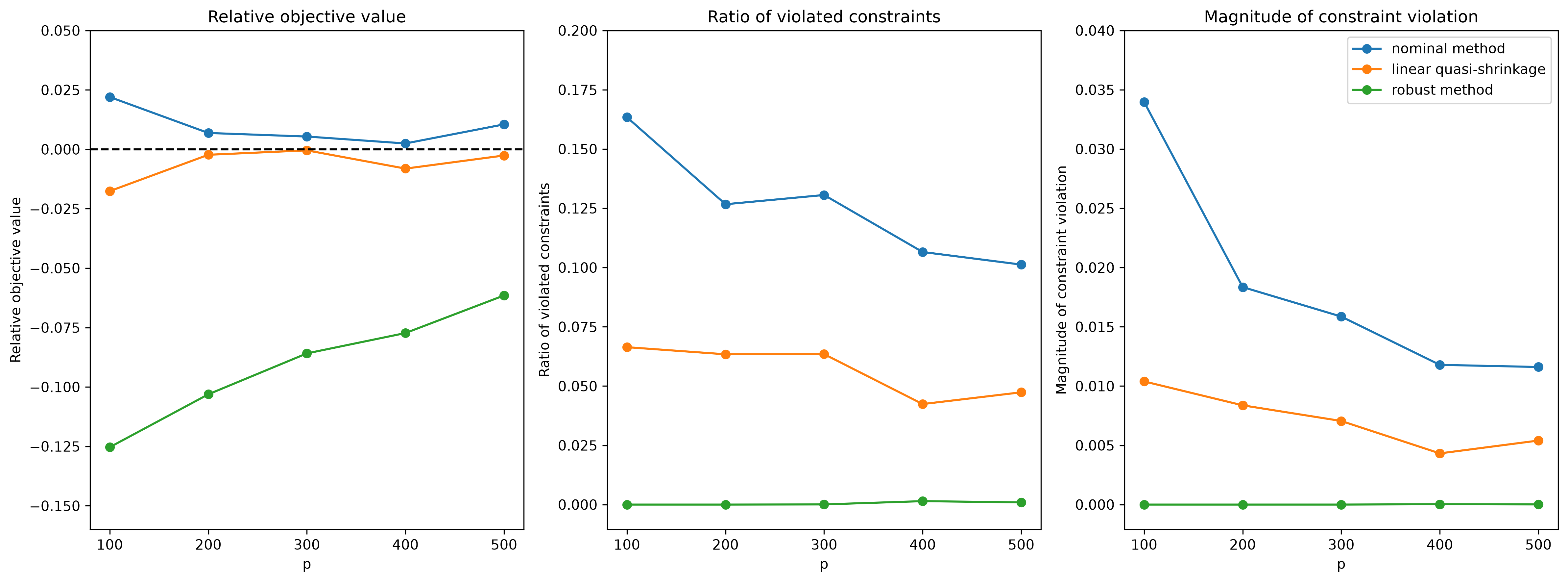}
\caption{Performance for $c=1,\, n =3,\, \sigma = 3$}
\label{c10_n3_s3}
\end{subfigure}

\end{figure}

\begin{figure}

\begin{subfigure}
\centering

\includegraphics[width=1.05\textwidth]{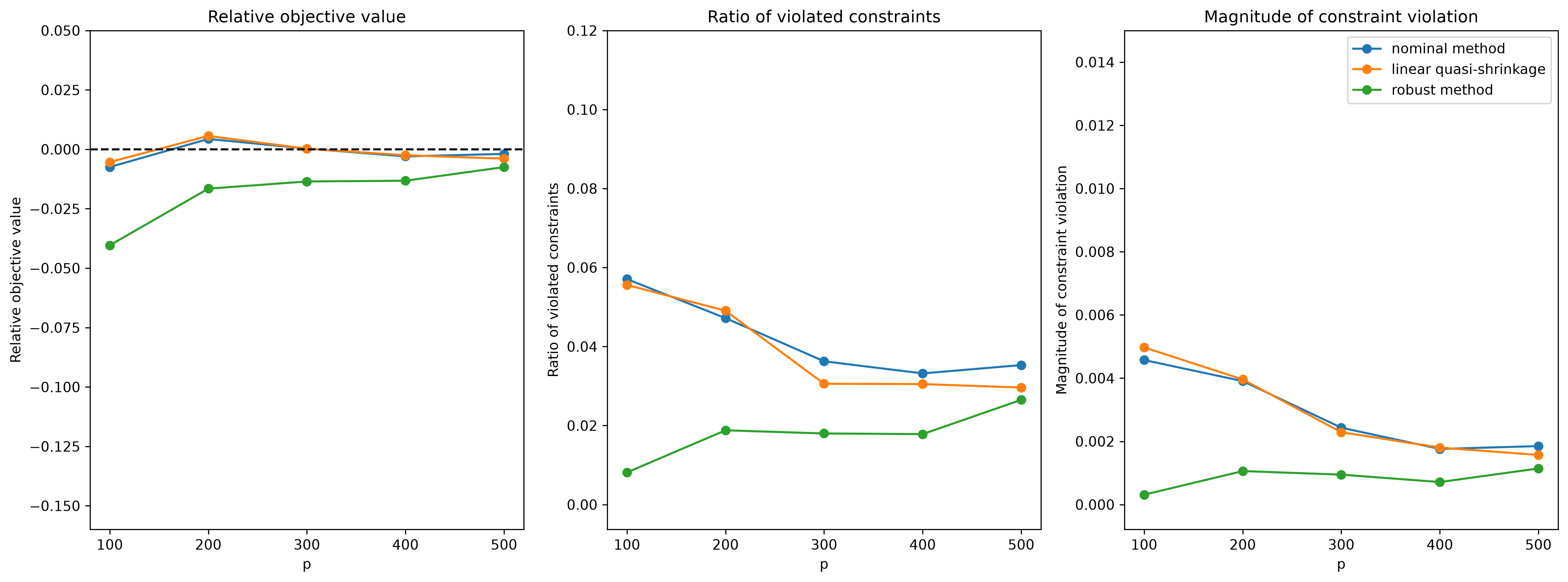}
\caption{Performance for  $c=1.5,\, n =3,\, \sigma =1$}
\label{c15_n3_s1}
\end{subfigure}
  
\begin{subfigure}
\centering

\includegraphics[width=1.05\textwidth]{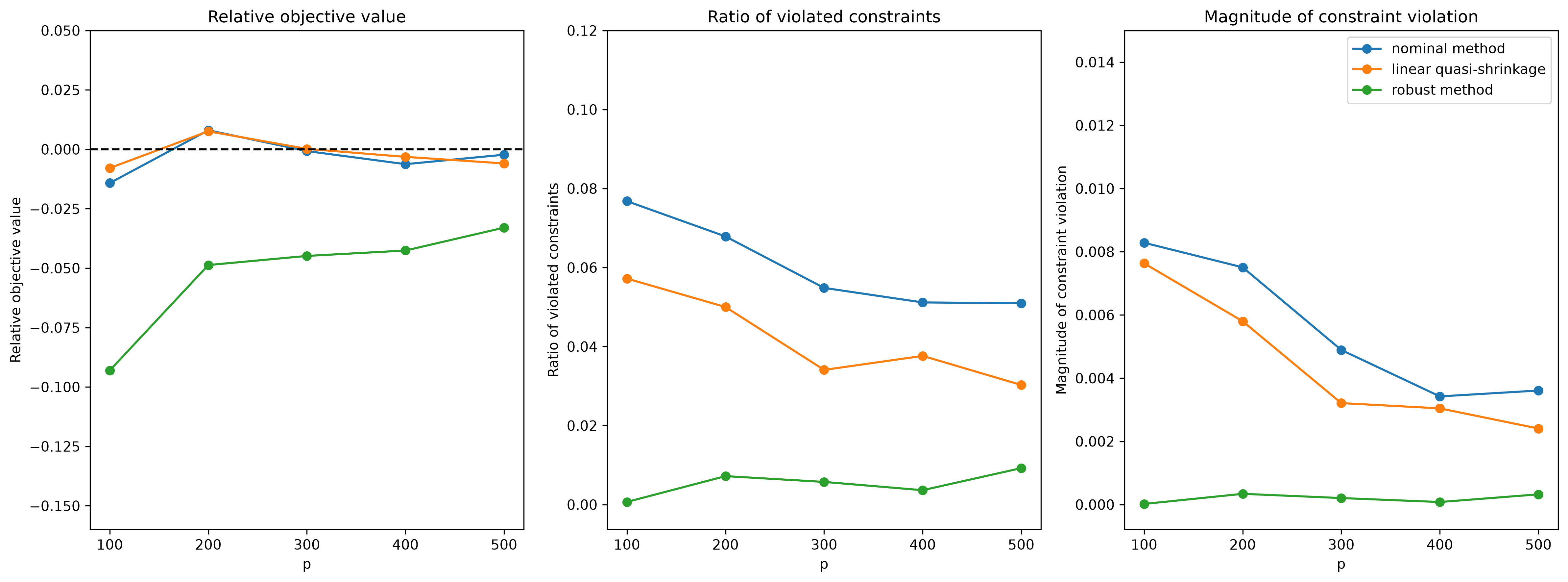}
\caption{Performance for  $c=1.5,\, n =3,\, \sigma =2$}
\label{c15_n3_s2}
\end{subfigure}

\begin{subfigure}
\centering

\includegraphics[width=1.05\textwidth]{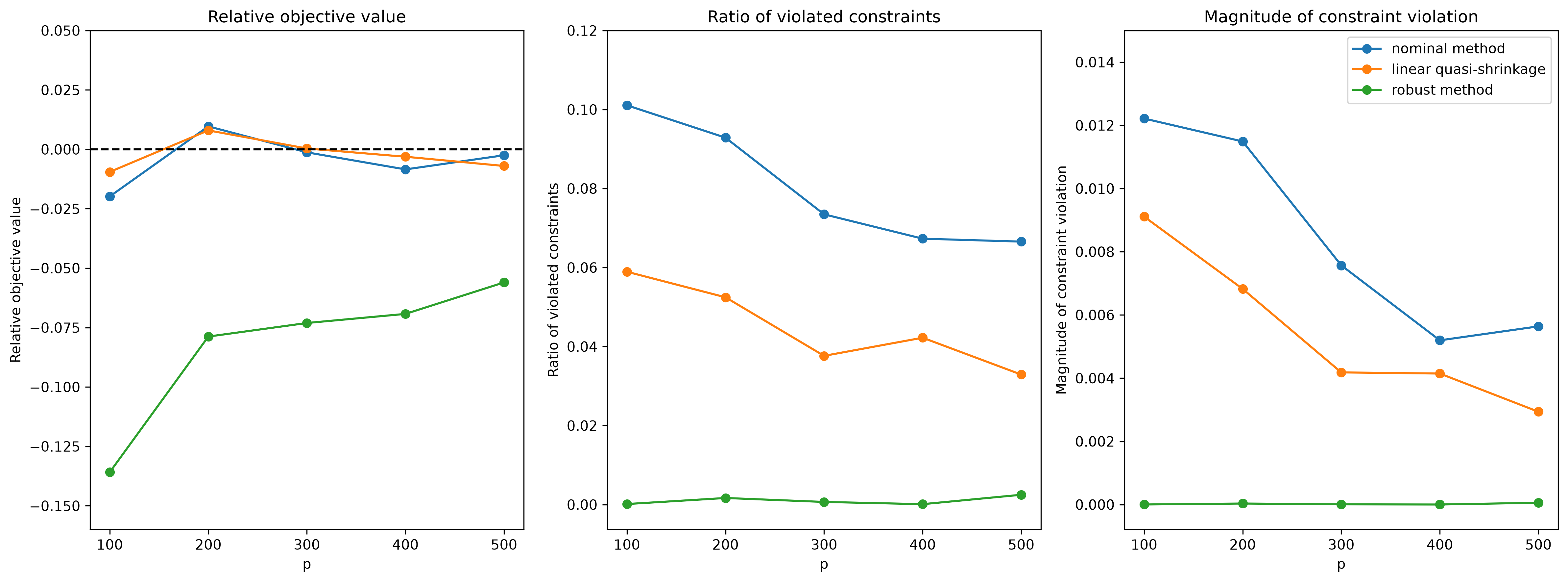}
\caption{Performance for  $c=1.5,\, n =3,\, \sigma =3$}
\label{c15_n3_s3}
\end{subfigure}
\end{figure}

Aforementioned patterns are also observed in Figures \ref{c10_n3_s1}, \ref{c10_n3_s2}, \ref{c10_n3_s3} with $c=1,\,n=3$ and $\sigma =1 , 2, 3$, respectively; and in Figures \ref{c15_n3_s1}, \ref{c15_n3_s2}, \ref{c15_n3_s3} with $c =1.5,\, n=3$ and $\sigma =1,2, 3$, respectively. Especially, for $c=1.5$, the nominal and linear quasi-shrinkage methods have
similar performance across the three metrics when $\sigma=1$, and the advantage of
linear quasi-shrinkage in constraint violation is more evident
at higher noise levels.

\begin{figure}
\centering

\includegraphics[width=1.05\textwidth]{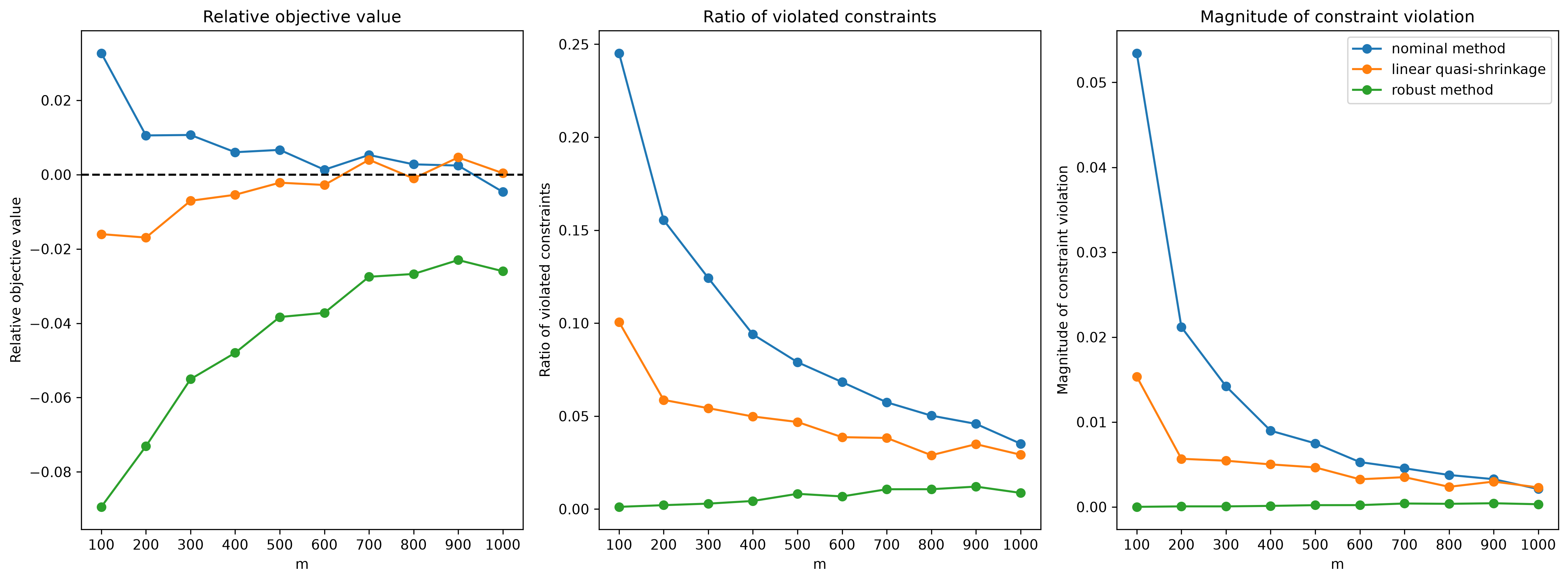}
\caption{Performance for  $p=500,\, n =3,\, \sigma =2$ and varying $m$}
\label{vary_m_p500_s2_n3}
\end{figure}

In Figure \ref{vary_m_p500_s2_n3}, we further investigate how the ratio $c =m/p$ influences the performance of the linear quasi-shrinkage method. We fix $p=500,\,n=3,\, \sigma =2$ and increase $m$ from $100$ to $1000$ (equivalently, $c = m/p$ ranges from $0.2$ to $2$) with steps of $100$. The results show that the linear quasi-shrinkage method achieves objective values close to the true objective value throughout the range considered, and its advantages over the nominal method in terms of constraint violation and over the robust method in terms of objective value are greater when the number of constraints $m$ is much smaller than the number of variables $p$, i.e., when $c =m/p$ is small. As $m$ and hence $c$ increases, both the ratio of violated constraints and the magnitude of constraint violation generally decrease for the nominal and linear quasi-shrinkage methods, and their performance becomes increasingly similar. The robust method maintains low constraint violations, while its relative objective value generally approaches zero but remains below the other two methods.

\begin{figure}

\begin{subfigure}
\centering
\includegraphics[width=1.05\textwidth]{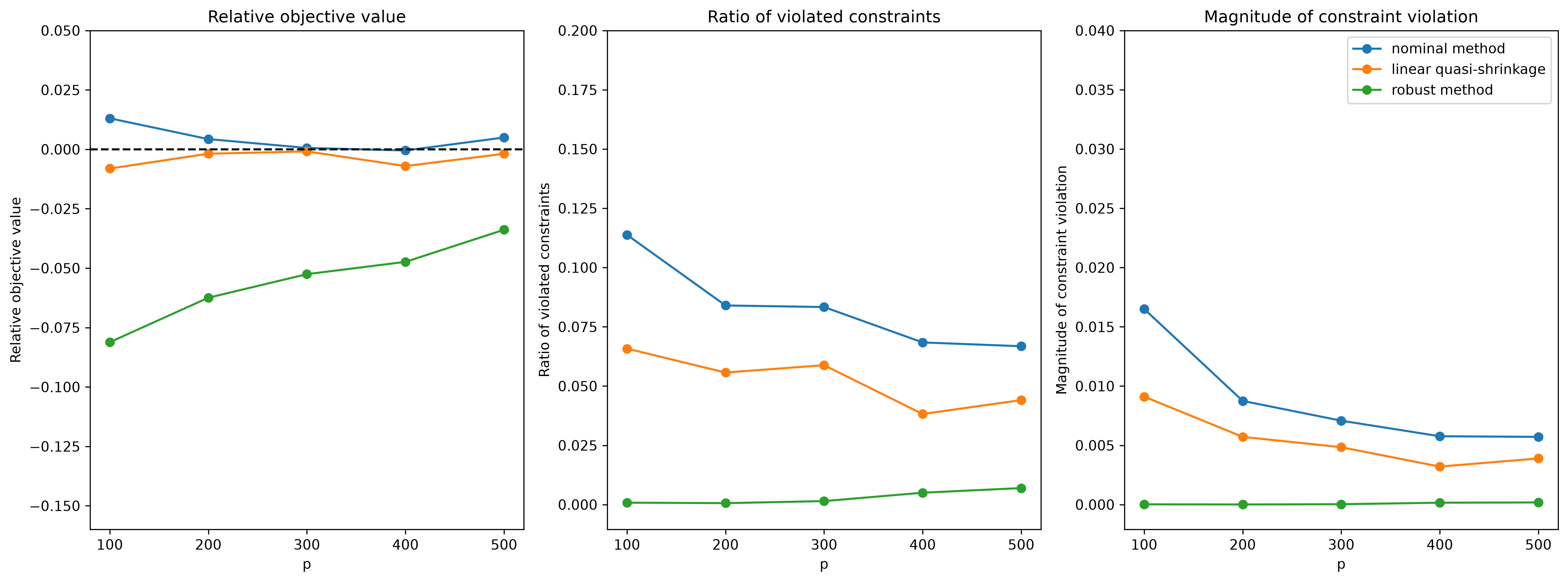}
\caption{Performance for $c=1, n=5, \sigma=2$}
\label{c10_n5_s2}
\end{subfigure}

\begin{subfigure}
\centering
\includegraphics[width=1.05\textwidth]{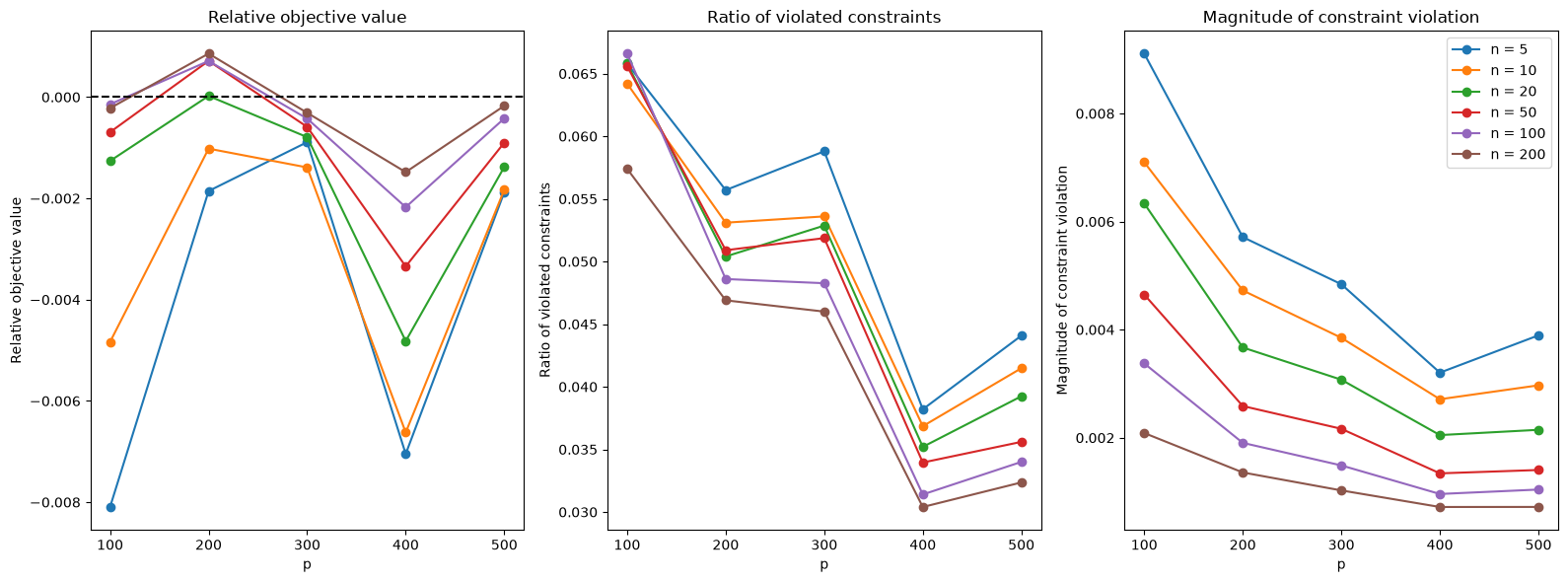}
\caption{Performance of linear quasi-shrinkage for $c=1,\, \sigma = 2$ and varying $n$}
\label{vary_n_c1_s2}
\end{subfigure}

\begin{subfigure}
\centering
\includegraphics[width=1.05\textwidth]{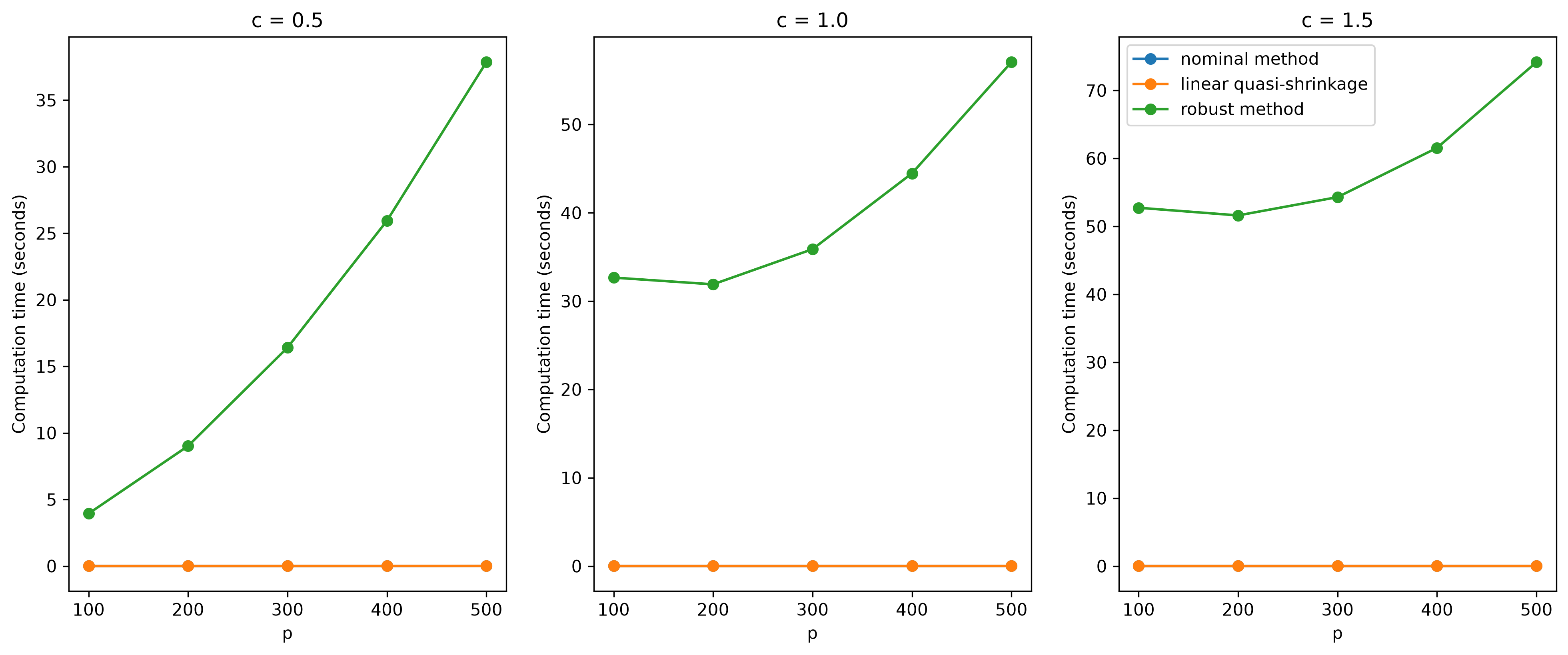}
\caption{Computation time for $n=3,\, \sigma = 2$}
\label{time_n3_s2}
\end{subfigure}

\end{figure}

To study the impact of the sample size $n$ on the performance, we present a representative case in Figure \ref{c10_n5_s2} where we take $c=1,\, n=5,\, \sigma = 2$, which exhibits a similar behavior to that of $n=3$. To facilitate further comparison, Figure~\ref{vary_n_c1_s2} shows results across a range of sample sizes for the linear quasi-shrinkage method with $c=1,\, \sigma = 2$. Although we do not include comparisons with the other methods here, one can expect that as $n$ increases, the performance of all three methods improves, and the linear quasi-shrinkage and nominal methods eventually converge in behavior. As shown in Figure \ref{vary_n_c1_s2}, increasing the sample size leads to more stable objective values that remain close to the true value. In contrast, the ratio of violated constraints and the magnitude of constraint violation decrease substantially. Overall, this further highlights the stability of the linear quasi-shrinkage method with respect to the sample size. 

Regarding the computation time, the robust method is expected to take significantly longer than both the nominal method and linear quasi-shrinkage method. This is because the latter two remain linear programming problems, while the robust method in this case results in a second-order cone programming problem and its reported computation time includes the cross-validation procedure for selecting the radius $\gamma$. We present Figure \ref{time_n3_s2} to illustrate this difference. Note that the nominal method and linear quasi-shrinkage method take approximately $10^{-3}$ to $10^{-2}$ seconds per instance, making the computation time negligible compared to the robust method and visually indistinguishable in the plots.

In conclusion, our linear quasi-shrinkage method stands out for its resilience against high levels of noise and varied sample sizes, consistently maintaining a more stable performance in objective value, the ratio of violated constraints and magnitude of constraint violation when compared to the nominal method and robust method. It should be noted that our linear quasi-shrinkage method performs well even for larger $c$, e.g., 2, although its advantage compared to the nominal method  is smaller in these settings. Furthermore, the linear quasi-shrinkage method preserves linearity of the problem, resulting in a low computational cost.


\section{Conclusion}\label{Sec:Conclusion}
In this paper, we investigate high-dimensional optimization problems featuring linear constraints under uncertainty, where the numbers of constraints and variables are both large and the parameter matrix is only known via noisy samples with small sample size. We assume additive noise for the way the true matrix is perturbed and consider column (row)-correlated noise matrices. To address such problems, we apply a linear quasi-shrinkage method to the parameter matrix. More precisely, our goal is to obtain an estimator of the true matrix, which takes the form of a linear combination of the matrix based on observations and the matrix representing prior information. Then, the estimator matrix is used as input to solve the optimization problem. The target matrix can incorporate available structural or prior information, while the shrinkage coefficients are consistently estimated from the observations without additional tuning or calibration. Under the stated assumptions, the asymptotic normalized Frobenius loss of the proposed estimator is no greater than that of the nominal estimator. Moreover, simulation results show that the linear quasi-shrinkage method exhibits a consistent stability in the face of increasing noise levels and different sample sizes, achieving a low ratio of violated constraints and magnitude of constraint violation while keeping the objective value close to the true objective value. 

A key challenge of our linear quasi-shrinkage method could lie in the selection of an appropriate target matrix, which influences the solution quality. Exploring strategies for selecting an appropriate target matrix, for example, by exploiting domain-specific knowledge, could provide a more robust foundation for our linear quasi-shrinkage method. 
Thus, future research could focus on selecting an appropriate target matrix for real-world optimization problems. Other future research directions could include incorporating the objective function, the right-hand side of the constraints, and the decision variables into the method, especially in settings involving discrete decision variables. Another direction is to explore different noise models on the parameter matrix, allowing for more general and realistic forms of uncertainty. 
\bibliography{ref}
\newpage
\newpage
\section*{Appendix}
\appendix
\renewcommand{\thesection}{\Alph{section}.\arabic{section}}
\begin{appendix}

\section{Supplemental proofs }\label{app:supplement proofs}

To assist with the proof of Theorem \ref{thm:cor-asymptotics}, we introduce the following lemma, which is a natural consequence of Borel-Cantelli lemma (see, e.g. \cite{shiryaev2016probability}).
\begin{lemma}\label{lemma:Borel Cantelli}
$X_m$ ($m=1,2,\dots$) is a sequence of random variables. If 
\begin{equation}
    \sum\limits_{m=1}^\infty\mathbb{P}\left(\left|X_m - \mathbb{E}\left(X_m\right)\right| > \frac1{m^\delta}\right) <\infty,
\end{equation}
    for some $\delta >0$, then $X_m - \mathbb{E}\left(X_m\right) \xrightarrow{a.s.} 0$.
\end{lemma}

\begin{lemma}\label{APP:correlated KSLLN of A, U with Sigma}
 $\frac{1}{mp}\left|\dotprod{\bS^{1/2}\bE}{\bA}\right| \xrightarrow{a.s} 0$ and $\frac{1}{mp}\left|\dotprod{\bS^{1/2}\bE}{\bU}\right| \xrightarrow{a.s} 0$ under Assumptions \eqref{ass:bounded trace}, \eqref{ass:high dimensional regime} and \eqref{ass:correlated noise model}.         
\end{lemma}
\begin{proof}

By Markov inequality, for any $\delta > 0$, 
\begin{equation}\label{eq: markov}
\sum\limits_{m=1}^{\infty}\mathbb{P}\left(\frac{1}{mp}\left|\dotprod{\bS^{1/2}\bE}{\bA}\right| > \frac1{(mp)^\delta} \right) 
\leq \sum\limits_{m=1}^{\infty}\frac{\mathbb{E}\left(\left|\dotprod{\bS^{1/2}\bE}{\bA}\right|^{4+\varepsilon}\right) }{(mp)^{4+\varepsilon-\delta(4+\varepsilon)}}
\end{equation}
where
\begin{equation}\label{eq: inner SEA}
  \dotprod{\bS^{1/2}\bE}{\bA} = \sum\limits_{i=1}^m \sum\limits_{j=1}^p  \bE^{ij}\left(\bS^{1/2}\bA\right)^{ij} ,
\end{equation}
$\bE^{ij}$ and   $\left(\bS^{1/2}\bA\right)^{ij}$ are the $i,j$-th entry of $\bE$ and $\left(\bS^{1/2}\bA\right)$, respectively.

Since $\bE$ is a matrix of i.i.d. entries with zero mean and unit variance by Assumption \eqref{ass:correlated noise model}, we have
\begin{equation}\label{eq: rosenthal}
    \mathbb{E}\left|\sum\limits_{i=1}^m \sum\limits_{j=1}^p \left( \bE^{ij}\left(\bS^{1/2}\bA\right)^{ij} \right)\right|^{4+\varepsilon} \leq  C\left(\left(\sum\limits_{i=1}^m \sum\limits_{j=1}^p\left(\left(\bS^{1/2}\bA\right)^{ij}\right)^2\right)^{2+\frac{\varepsilon}{2}} + \sum\limits_{i=1}^m \sum\limits_{j=1}^p\left|\left(\bS^{1/2}\bA\right)^{ij}\right|^{4+\varepsilon}\right)
\end{equation}
by applying Rosenthal's inequality (see, e.g., \cite{rosenthal1970subspaces}), where $C$ is a constant depending on $4+\varepsilon$ and $\mathbb{E}\left(\left|\bE^{ij}\right|^{4+\varepsilon}\right)$ which is uniformly bounded over $m$ and $p$ according to Assumption \eqref{ass:correlated noise model}. Note that
\begin{equation}\label{eq: rosenthal 1}
  \sum\limits_{i=1}^m \sum\limits_{j=1}^p\left(\left(\bS^{1/2}\bA\right)^{ij}\right)^2= \mathrm{tr}\left(\bA^\top\bS\bA\right), 
\end{equation}
and
\begin{equation}\label{eq: rosenthal 2}
  \sum\limits_{i=1}^m \sum\limits_{j=1}^p\left|\left(\bS^{1/2}\bA\right)^{ij}\right|^{4+\varepsilon} \leq  \left(\sum\limits_{i=1}^m \sum\limits_{j=1}^p\left|\left(\bS^{1/2}\bA\right)^{ij}\right|^2\right) ^{2+\frac{\varepsilon}{2}}  = \left(\mathrm{tr}\left(\bA^\top\bS\bA\right)\right)^{2+\frac{\varepsilon}{2}}.
\end{equation}
Combining  \eqref{eq: inner SEA}, \eqref{eq: rosenthal}, \eqref{eq: rosenthal 1} and \eqref{eq: rosenthal 2}, we obtain
\begin{equation}
     \mathbb{E}\left(\left|\dotprod{\bS^{1/2}\bE}{\bA}\right|^{4+\varepsilon}\right) = O\left(\left(\mathrm{tr}\left(\bA^\top\bS\bA\right)\right)^{2+\frac{\varepsilon}{2}}\right) = O\left(\left(\mathrm{tr}\left(\bA^\top\bA\right)\mathrm{tr}\left(\bS\right)\right)^{2+\frac{\varepsilon}{2}}\right)
\end{equation}
where the second bound follows from the trace inequality. Consequently, 
\begin{equation}
\frac{\mathbb{E}\left(\left|\dotprod{\bS^{1/2}\bE}{\bA}\right|^{4+\varepsilon}\right) }{(mp)^{4+\varepsilon-\delta(4+\varepsilon)}}  = O\left(\left(\frac{\mathrm{tr}\left(\bA^\top\bA\right)}{mp}\right)^{2+\frac{\varepsilon}{2}} \left(\frac{\mathrm{tr}\left(\bS\right)}{m}\right)^{2+\frac{\varepsilon}{2}} \frac{(mp)^{\delta(4+\varepsilon)}}{p^{2+\frac{\varepsilon}{2}}} \right). 
\end{equation}
Under Assumptions \eqref{ass:bounded trace}, \eqref{ass:high dimensional regime} and \eqref{ass:correlated noise model}, we have respectively that
\begin{equation}
    \quad \sup\limits_{m} \frac{\mathrm{tr}\left(\bS\right)}m < \infty, \quad \frac{m}p \rightarrow c \in [0, \infty) \quad \text{as}\quad m, p\rightarrow \infty   \quad \text{and} \quad\sup\limits_{m}\frac1{mp} \mathrm{tr}\left(\bA^\top\bA\right) < \infty.
\end{equation}
Thus, recalling \eqref{eq: markov}, we obtain
\begin{equation}
\sum\limits_{m=1}^{\infty}\mathbb{P}\left(\frac{1}{mp}\left|\dotprod{\bS^{1/2}\bE}{\bA}\right| > \frac1{(mp)^\delta} \right)  \leq \sum\limits_{m=1}^{\infty}\frac{\mathbb{E}\left(\left|\dotprod{\bS^{1/2}\bE}{\bA}\right|^{4+\varepsilon}\right) }{(mp)^{4+\varepsilon-\delta(4+\varepsilon)}},
\end{equation}
where 
\[ \frac{\mathbb{E}\left(\left|\dotprod{\bS^{1/2}\bE}{\bA}\right|^{4+\varepsilon}\right) }{(mp)^{4+\varepsilon-\delta(4+\varepsilon)}} = O\left(\frac1{m^{2+\frac{\varepsilon}2 -2\delta(4+\varepsilon)}}\right). \]
Hence,
\begin{equation}
    \sum\limits_{m=1}^{\infty}\mathbb{P}\left(\frac{1}{mp}\left|\dotprod{\bS^{1/2}\bE}{\bA}\right| > \frac1{(mp)^\delta} \right)  <\infty
\end{equation}
for $0 < \delta < \frac14\left(1-\frac{2}{4+\varepsilon}\right)$. Therefore, by Lemma \ref{lemma:Borel Cantelli},
\begin{equation}
 \frac{1}{mp}\left|\dotprod{\bS^{1/2}\bE}{\bA}\right| \xrightarrow{a.s} 0 
\end{equation}
as $m, p\rightarrow \infty$ with $m/p \rightarrow c \in [0, \infty)$. Similarly,
\begin{equation}
 \frac{1}{mp}\left|\dotprod{\bS^{1/2}\bE}{\bU}\right| \xrightarrow{a.s} 0  
\end{equation}
as $m, p\rightarrow \infty$ with $m/p \rightarrow c \in [0, \infty)$.

\end{proof}

\section{Proof of Theorem \ref{theorem: asymptotic loss}  }\label{App: proof of asymptotic loss}
\begin{proof}

Recalling \eqref{eq: bounded denominator} and applying the asymptotic results in Theorem \ref{thm:cor-asymptotics}, we have
\begin{equation*}
\begin{aligned}
    \frac1{mp}L^*_{m}  \overset{a.e.}{\sim} & \quad- \frac1{mp}\left[\frac{\left(\dotprod{\bA}{\bA}\dotprod{\bU}{\bU} - \dotprod{\bA}{\bU}^2\right)^2}{\dotprod{\bA}{ \bA}\dotprod{\bU}{\bU}
     + \frac{p}n \dotprod{\bd\Sigma_m}{\mathbf{I}_m} \dotprod{\bU}{\bU}- \dotprod{\bA}{\bU}^2}\right]\frac1{\dotprod{\bU}{\bU}} \\[1ex]
    &\quad - \frac1{mp}\left[\frac{ \dotprod{\bA}{\bU}^2 }{ \dotprod{\bU}{\bU} } - \dotprod{\bA}{\bA}
\right]  \\[1ex]
     & =  - \frac1{mp}\left[\frac{\left(\dotprod{\bA}{\bA}\dotprod{\bU}{\bU}+ \frac{p}n \dotprod{\bd\Sigma_m}{\mathbf{I}_m} \dotprod{\bU}{\bU}- \dotprod{\bA}{\bU}^2 - \frac{p}n \dotprod{\bd\Sigma_m}{\mathbf{I}_m}\dotprod{\bU}{\bU}\right )^2 }{\dotprod{\bA}{ \bA}\dotprod{\bU}{\bU} + \frac{p}n \dotprod{\bd\Sigma_m}{\mathbf{I}_m}\dotprod{\bU}{\bU} - \dotprod{\bA}{\bU}^2}\right]\\[1ex]
       &\quad \times\frac1{\dotprod{\bU}{\bU}}
   - \frac1{mp} \left[\frac{ \dotprod{\bA}{\bU}^2 }{ \dotprod{\bU}{\bU} } - \dotprod{\bA}{\bA}
\right] \\[1ex]
     & = -\frac1{mp}\left[\dotprod{\bA}{ \bA}\dotprod{\bU}{\bU} + \frac{p}n \dotprod{\bd\Sigma_m}{\mathbf{I}_m}\dotprod{\bU}{\bU} - \dotprod{\bA}{\bU}^2 - 2\,\frac{p}n \dotprod{\bd\Sigma_m}{\mathbf{I}_m}\dotprod{\bU}{\bU}\right.\\[1ex]
     &\quad \left. + \frac{\frac{p^2}{n^2} \dotprod{\bd\Sigma_m}{\mathbf{I}_m}^2\dotprod{\bU}{\bU}^2}{\dotprod{\bA}{ \bA}\dotprod{\bU}{\bU} + \frac{p}n \dotprod{\bd\Sigma_m}{\mathbf{I}_m}\dotprod{\bU}{\bU} - \dotprod{\bA}{\bU}^2} \right]\frac1{\dotprod{\bU}{\bU}}\\[1ex]
     &\quad - \frac1{mp} \left[\frac{ \dotprod{\bA}{\bU}^2 }{ \dotprod{\bU}{\bU} } - \dotprod{\bA}{\bA}
\right] \\[1ex]
     & =\frac1{mp}\frac{\left(\dotprod{\bA}{ \bA}\dotprod{\bU}{\bU}-\dotprod{\bA}{\bU}^2\right)\frac{p}n \dotprod{\bd\Sigma_m}{\mathbf{I}_m}\dotprod{\bU}{\bU} }{\dotprod{\bA}{ \bA}\dotprod{\bU}{\bU} + \frac{p}n \dotprod{\bd\Sigma_m}{\mathbf{I}_m}\dotprod{\bU}{\bU} - \dotprod{\bA}{\bU}^2} \frac{1 }{\dotprod{\bU}{\bU}}\\[1ex]
     & = \frac{\dotprod{\bA}{ \bA}\dotprod{\bU}{\bU}-\dotprod{\bA}{\bU}^2}{\dotprod{\bA}{ \bA}\dotprod{\bU}{\bU} + \frac{p}n \dotprod{\bd\Sigma_m}{\mathbf{I}_m}\dotprod{\bU}{\bU} - \dotprod{\bA}{\bU}^2}\frac{1}{nm} \dotprod{\bd\Sigma_m}{\mathbf{I}_m}\\[1ex]
     & = \frac{\alpha^*}{nm}\dotprod{\bd\Sigma_m}{\mathbf{I}_m}
\end{aligned}
\end{equation*} \qedhere
\end{proof}

\section{Some definitions and theorems}\label{App: Miscellaneous}


\begin{definition}[Consistent estimator]\label{def: consistent estimator}
An estimator \(\hat{\theta}_m\) is said to be consistent for a parameter \(\theta_m\) as \(m\) increases, if the estimator converges almost surely to the true value of the parameter. Formally:

\[
\mathbb{P}\left(\lim\limits_{m\to\infty}\left(\hat{\theta}_m -\theta_m\right)=0\right)=1 \quad \overset{def.}{\Leftrightarrow}\quad \hat{\theta}_m - \theta_m \xrightarrow{a.s.} 0 \quad \text{as} \quad m \to \infty.
\]

\end{definition}

\vspace{2em}

We give a brief proof of the following result, which is a special case of Theorem 4.2 in \cite{mohajerin2018data}. 

\begin{theorem}\label{theorem: wasserstein}
 Define $\mathcal{B}(\hat{\mathbb{P}}, \gamma) = \left\{\mathbb{P}: \|\mathbb{P} -  \hat{\mathbb{P}}\|_W \leq \gamma \right\}$ where   $\hat{\mathbb{P}}$ is the empirical distribution over $\mathbf{a}$  and $\|\cdot\|_W$ is the 1-Wasserstein distance, defined as
\begin{equation*}
    d_{W}(\mathbb{P}_1, \mathbb{P}_2) = \sup\limits_{f \in \mathcal{L}} \left\{\int\limits_{\Xi} f(\boldsymbol\xi)\mathbb{P
    }_1(d\boldsymbol\xi) - \int\limits_{\Xi} f(\boldsymbol\xi)\mathbb{P}_2(d\boldsymbol\xi)   \right\}
\end{equation*}
where $\Xi=\mathbb{R}^p$ and $\mathcal{L}$ denotes the set of Lipschitz functions on $\Xi$ satisfying $|f(\boldsymbol\xi) - f(\boldsymbol\xi')|\leq \|\boldsymbol\xi - \boldsymbol\xi'\|$ for all $\boldsymbol\xi$, $\boldsymbol\xi'\in \Xi$ with respect to Euclidean norm. Then, we have
  \begin{equation*}
  \sup\limits_{\mathbb{P}\in \mathcal{B}(\hat{\mathbb{P}}, \gamma)}\left\{\mathbb{E}_{\mathbb{P}}\left(\mathbf{a}^\top\mathbf{x}\right)\right\} = \mathbb{E}_{\hat{\mathbb{P}}} (\mathbf{a}^\top \bx)  + \gamma \|\mathbf{x}\|. 
\end{equation*} 
    
\end{theorem}
\begin{proof}
According to the proof of Theorem 4.2 in \cite{mohajerin2018data}, for a linear function $l(\mathbf{a}) = \mathbf{a}^\top\mathbf{x}$, the worst-case expectation over a Wasserstein ball is given by:

\begin{equation}\label{theorem 4.2 formula}
        \sup\limits_{\mathbb{P}\in \mathcal{B}(\hat{\mathbb{P}}, \gamma)} \mathbb{E}_{\mathbb{P}} \left(\mathbf{a}^\top\mathbf{x}\right) = \inf\limits_{\lambda\geq 0} \left\{\lambda \gamma + \frac1n \sum\limits_{k=1}^n \sup\limits_{\mathbf{a}\in \mathbf{R}^p}\left\{\mathbf{a}^\top\mathbf{x} - \lambda \|\mathbf{a} - \tilde{\mathbf{a}}^k\|\right\}\right\}
    \end{equation}
where $\|\cdot\|$ is the 2-norm and $\tilde{\mathbf{a}}^k$ denotes the $k$-th observation. Using H\"older's inequality, for any $k\in \{1,2, \dots, n\}$,
\begin{equation*}
    |\mathbf{a}^\top\mathbf{x} - 
   (\tilde{\mathbf{a}}^k)^\top\mathbf{x}| \leq \|\mathbf{a} -  \tilde{\mathbf{a}}^k \|\cdot \|\mathbf{x}\|.
\end{equation*}
Hence, 
\begin{equation}\label{holder inequality}
\mathbf{a}^\top\mathbf{x} \leq 
   (\tilde{\mathbf{a}}^k)^\top\mathbf{x} + \|\mathbf{a} -  \tilde{\mathbf{a}}^k \|\cdot \|\mathbf{x}\|.
\end{equation} 
Applying \eqref{holder inequality}, we have
\begin{equation}\label{upper bound}
    \mathbf{a}^\top\mathbf{x} - \lambda \|\mathbf{a} -  \tilde{\mathbf{a}}^k\|\leq  (\tilde{\mathbf{a}}^k)^\top\mathbf{x} + \left(\|\mathbf{x}\| - \lambda\right)\cdot\|\mathbf{a} -  \tilde{\mathbf{a}}^k\|.
\end{equation}
The right-hand side of \eqref{upper bound} is maximized at $\mathbf{a} = \tilde{\mathbf{a}}_k $ when   $  \lambda \geq \|\mathbf{x}\| $, i.e.,
\begin{equation*}
    \sup\limits_{\mathbf{a}\in \mathbf{R}^p}\left(\mathbf{a}^\top\mathbf{x} - \lambda \|\mathbf{a} - \tilde{\mathbf{a}}_k\|\right) = (\tilde{\mathbf{a}}^k)^\top\mathbf{x}.
\end{equation*}
Thus, for $\lambda \geq \|\mathbf{x}\|$,  the right-hand side of \eqref{theorem 4.2 formula} is minimized when $\lambda = \|\mathbf{x}\|$, i.e.,  
\begin{equation*}
    \inf\limits_{\lambda\geq \|\mathbf{x}\| } \left\{\lambda \gamma + \frac1n\sum\limits_{k=1}^n\sup\limits_{\mathbf{a}\in \mathbf{R}^p}\left(\mathbf{a}^\top\mathbf{x} - \lambda \|\mathbf{a} - \tilde{\mathbf{a}}_k\|\right)\right\} = \frac1n\sum\limits_{k=1}^n (\tilde{\mathbf{a}}^k)^\top\mathbf{x} + \gamma\|\mathbf{x}\|.
 \end{equation*}
When $\lambda <\|\mathbf{x}\|$ and $\|\bx\| \neq 0$, let $\mathbf{a}_t = \tilde{\mathbf{a}}^k + t\frac{\bx}{\|\bx\|}$. Then, 
\[\mathbf{a}_t^\top \bx - \lambda\|\mathbf{a}_t - \tilde{\mathbf{a}}^k\| = \left(\tilde{\mathbf{a}}^k\right)^\top \bx + t\left(\|\bx\| - \lambda\right) \rightarrow \infty\]
as $t\rightarrow \infty$. Therefore, 
\begin{equation*}
     \sup\limits_{\mathbf{a}\in \mathbb{R}^p} \left\{\mathbf{a}^\top\mathbf{x} - \lambda \|\mathbf{a} -  (\tilde{\mathbf{a}}^k)\|\right\} = \infty
\end{equation*}
according to \eqref{upper bound}. Therefore, we conclude that
\begin{equation*}
    \inf\limits_{\lambda \geq 0}  \left\{\lambda \gamma + \frac1n\sum\limits_{k=1}^n\sup\limits_{\mathbf{a}\in \mathbf{R}^p}\left\{\mathbf{a}^\top\mathbf{x} - \lambda \|\mathbf{a} - \tilde{\mathbf{a}}_k\|\right\}\right\} = \frac1n\sum\limits_{k=1}^n(\tilde{\mathbf{a}}^k)^\top\mathbf{x} +\gamma \|\mathbf{x}\|.
 \end{equation*}
Hence, according to \eqref{theorem 4.2 formula},
\begin{equation*}
     \sup\limits_{\mathbb{P}\in \mathcal{B}(\hat{\mathbb{P}}, \gamma)} \mathbb{E}_{\mathbb{P}} \left(\mathbf{a}^\top\mathbf{x}\right) = \frac1n\sum\limits_{k=1}^n (\tilde{\mathbf{a}}^k)^\top\mathbf{x} + \gamma \|\mathbf{x}\|.\qedhere
\end{equation*}

\end{proof}

\section{Solve DRO using cross-validation}\label{App:cross-validation}

This section describes the cross-validation procedure used to select the robustness radius $\gamma$ for the DRO method and outlines the algorithm.

\begin{enumerate}

    \item \textbf{Construct the initial candidate set for $\gamma$.}

    Define the initial set of multipliers as
    \[
    \Lambda_0
    =
    \{0,0.1,0.2,\ldots,2.0\}.
    \]
    The corresponding candidate set for $\gamma$ is
    \[
    \Gamma_0
    =
    \left\{
    \lambda\widehat{\sigma}_m:\lambda\in\Lambda_0
    \right\}.
    \]
where $\widehat{\sigma}_m$ is the estimated standard deviation of the noise. 
    \item \textbf{Construct the leave-one-out training samples.}

    For each validation index $j\in\{1,\ldots,n\}$, leave out the observation
    $\tilde{\mathbf{A}}_m^j$ and compute the training-sample average
    \[
    \bar{\mathbf{A}}_m^{(-j)}
    =
    \frac{1}{n-1}
    \sum_{\substack{k=1\\k\neq j}}^n
    \tilde{\mathbf{A}}_m^k
    =
    \frac{
    \sum\limits_{k=1}^n\tilde{\mathbf{A}}_m^k
    -
    \tilde{\mathbf{A}}_m^j
    }{n-1}.
    \]

    \item \textbf{Solve the robust problem on each training sample.}

    For every $\gamma\in\Gamma_0$ and every validation index $j$, solve
    \[
    \mathbf{x}_{\gamma}^{(-j)}
    \in
    \arg\max_{\mathbf{x}\geq\mathbf{0}}
    \left\{\mathbf{c}^{\top}\mathbf{x}: \quad\bar{\mathbf{A}}_m^{(-j)}\mathbf{x}
    +
    \gamma\|\mathbf{x}\|_2\mathbf{1}_m
    \leq
    \mathbf{b}\right\}.    
    \]
    If this problem is not solved successfully for a particular candidate
    $\gamma$ in any fold, that candidate is excluded from consideration.

    \item \textbf{Evaluate the solution on the validation observation}

    For each constraint $i\in\{1,\ldots,m\}$, define the validation residual
    \[
    r_{ij}(\gamma)
    =
    \left(
    \widetilde{\mathbf{A}}_m^j
    \mathbf{x}_{\gamma}^{(-j)}
    -
    \mathbf{b}
    \right)_i.
    \]
    A validation constraint is considered violated if
    \[
    r_{ij}(\gamma)>\varepsilon_{\mathrm{feas}},
    \]
    where $\varepsilon_{\mathrm{feas}}=10^{-6}$ is the numerical feasibility
    tolerance. The ratio of violated constraints in validation fold $j$ is
    \[
    v_j(\gamma)
    =
    \frac{1}{m}
    \sum_{i=1}^m
    \mathbbm{1}
    \left\{
    r_{ij}(\gamma)>\varepsilon_{\mathrm{feas}}
    \right\}.
    \]
    The corresponding validation objective value is
    \[
    q_j(\gamma)
    =
    \mathbf{c}^{\top}\mathbf{x}_{\gamma}^{(-j)}.
    \]

    \item \textbf{Aggregate the validation results.}

    For each candidate $\gamma$, calculate the average ratio of violated constraints for all folds
    \[
    \bar{v}(\gamma)
    =
    \frac{1}{n}
    \sum_{j=1}^n v_j(\gamma),
    \]
    and the average objective value
    \[
    \bar{q}(\gamma)
    =
    \frac{1}{n}
    \sum_{j=1}^n q_j(\gamma).
    \]
    Given a prespecified threshold for the ratio of violated constraints $\tau$, a candidate is
    considered eligible if
    \[
    \bar{v}(\gamma)\leq\tau.
    \]
    In the simulation, we use $\tau=0.05$.

    \item \textbf{Extend the candidate set when necessary.}

    First, the preceding steps are performed over the initial candidate set
    $\Gamma_0$. If no candidate in $\Gamma_0$ is eligible, or if the candidate
    selected from $\Gamma_0$ is its largest element $2\widehat{\sigma}$, an
    additional multiplier set
    \[
    \Lambda_1
    =
    \{2.1,2.2,\ldots,3.0\}
    \]
    is evaluated. The corresponding additional candidate set is
    \[
    \Gamma_1
    =
    \left\{
    \lambda\widehat{\sigma}_m:\lambda\in\Lambda_1
    \right\}.
    \]
    In this case, the final selection is performed over all eligible candidates
    in $\Gamma_0\cup\Gamma_1$. If the grid is not extended, the selection is
    performed over $\Gamma_0$ only.

    \item \textbf{Select the robustness radius.}

    Let $\Gamma_{\mathrm{elig}}$ denote the set of candidates that are solved
    successfully in every fold and satisfy the requirement for the average ratio of violated constraints. The cross-validated robustness radius is selected according to
    \[
    \widehat{\gamma}_{\mathrm{CV}}
    =
    \min\left(
    \arg\max_{\gamma\in\Gamma_{\mathrm{elig}}}
    \bar{q}(\gamma)
    \right).
    \]
    Thus, among the candidates whose average ratio of violated constraints does not
    exceed $\tau$, we select the candidate with the largest average 
    objective value. If several candidates have the same average objective
    value, the smallest $\gamma$ is selected.


    \item \textbf{Solve the final model using all observations.}

    Cross-validation is used only to select $\widehat{\gamma}_{\mathrm{CV}}$.
    After its selection, the final decision is obtained by solving the robust
    model once using the full-sample average:
    \[
    \widehat{\mathbf{x}}_{\mathrm{CV}}
    \in
    \arg\max_{\mathbf{x}\geq\mathbf{0}}
    \quad
    \left\{\mathbf{c}^{\top}\mathbf{x}, \quad 
    \bar{\mathbf{A}}_m\mathbf{x}
    +
    \widehat{\gamma}_{\mathrm{CV}}
    \|\mathbf{x}\|_2\mathbf{1}_m
    \leq
    \mathbf{b}\right\}.
    \]

\end{enumerate}

\section{Multiple targets}\label{App: multiple targets}
Consider incorporating additional target matrices $\bU^i$ $(i = 1,2\, \dots, l)$ in the formulation: 
\begin{equation}
    \mathbf{A}_m^* = \alpha \bbA + \sum\limits_{i=1}^l\beta_i \bU^i. 
\end{equation}
We continue using the Frobenius Loss:
\begin{equation}
   \min\limits_{\alpha, \beta_1, \dots, \beta_l}\left\| \mathbf{A}_m^* - \mathbf{A}\right\|^2.
\end{equation}
By taking the derivatives of $\left\| \mathbf{A}_m^* - \mathbf{A}\right\|^2$ with respect to $\alpha$ and each $\beta_i~(i=1, \dots, l)$ and setting them equal to zero, we arrive at the following linear system:

\begin{equation}
\begin{bmatrix}
\dotprod{\bbA}{\bbA} & \dotprod{\bbA}{\bU^1} & \cdots & \dotprod{\bbA}{\bU^l}\\
\dotprod{\bbA}{\bU^1} & \dotprod{\bU^1}{\bU^1} & \cdots & \dotprod{\bU^1}{\bU^l}\\
\vdots & \vdots & \ddots & \vdots \\
\dotprod{\bbA}{\bU^l}& \dotprod{\bU^l}{\bU^1} & \cdots & \dotprod{\bU^l}{\bU^l}
\end{bmatrix}
\begin{bmatrix}
\alpha \\
\beta_1\\
\vdots \\
\beta_l
\end{bmatrix}
=
\begin{bmatrix}
\dotprod{ \bbA}{\bA}\\
\dotprod{\bA}{\bU^1} \\
\vdots \\
\dotprod{ \bA}{\bU^l}
\end{bmatrix}.
\end{equation}
Under the assumption that $\bbA$ and $\mathbf{U}_m^i$ ($i =1,2, \dots, l$) are jointly linearly independent, the matrix on the left-hand side is invertible. Thus, the optimal weights are given by

\begin{equation}
    \begin{bmatrix}
(\alpha)_m^* \\
(\beta_1)_m^*\\
\vdots \\
(\beta_l)_m^*
\end{bmatrix} = \begin{bmatrix}
\dotprod{\bbA}{\bbA} & \dotprod{\bbA}{\bU^1} & \cdots & \dotprod{\bbA}{\bU^l}\\
\dotprod{\bbA}{\bU^1} & \dotprod{\bU^1}{\bU^1} & \cdots & \dotprod{\bU^1}{\bU^l}\\
\vdots & \vdots & \ddots & \vdots \\
\dotprod{\bbA}{\bU^l}& \dotprod{\bU^l}{\bU^1} & \cdots & \dotprod{\bU^l}{\bU^l}
\end{bmatrix}^{-1} 
\begin{bmatrix}
\dotprod{ \bbA}{\bA}\\
\dotprod{\bA}{\bU^1} \\
\vdots \\
\dotprod{ \bA}{\bU^l}
\end{bmatrix}.
\end{equation}
Following the same procedure as in Theorem \ref{thm:cor-asymptotics}, under analogous assumptions, we obtain `asymptotic oracle' weights: 
\begin{equation}
    \begin{bmatrix}
\alpha^* \\
\beta_1^*\\
\vdots \\
\beta_l^*
\end{bmatrix} = \begin{bmatrix}
\dotprod{\bA}{\bA} + \frac{p}{n}\dotprod{\bd\Sigma_m}{\mathbf{I}_m} & \dotprod{\bA}{\bU^1}  & \cdots & \dotprod{\bA}{\bU^l}\\
\dotprod{\bA}{\bU^1} & \dotprod{\bU^1}{\bU^1} & \cdots & \dotprod{\bU^1}{\bU^l}\\
\vdots & \vdots & \ddots & \vdots \\
\dotprod{\bA}{\bU^l}& \dotprod{\bU^l}{\bU^1} & \cdots & \dotprod{\bU^l}{\bU^l}
\end{bmatrix}^{-1} 
\begin{bmatrix}
\dotprod{ \bA}{\bA}\\
\dotprod{\bA}{\bU^1} \\
\vdots \\
\dotprod{ \bA}{\bU^l}
\end{bmatrix}.
\end{equation}
Consequently, the `bona fide' weights are given by 
\begin{equation}
    \begin{bmatrix}
\widehat{\alpha}^* \\
\widehat{\beta}_1^*\\
\vdots \\
\widehat{\beta}_l^*
\end{bmatrix} = \begin{bmatrix}
\dotprod{\bbA}{\bbA} & \dotprod{\bbA}{\bU^1} & \cdots & \dotprod{\bbA}{\bU^l}\\
\dotprod{\bbA}{\bU^1} & \dotprod{\bU^1}{\bU^1} & \cdots & \dotprod{\bU^1}{\bU^l}\\
\vdots & \vdots & \ddots & \vdots \\
\dotprod{\bbA}{\bU^l}& \dotprod{\bU^l}{\bU^1} & \cdots & \dotprod{\bU^l}{\bU^l}
\end{bmatrix}^{-1} 
\begin{bmatrix}
\dotprod{\bbA}{\bbA} - \frac{mp}{n}\widehat{\sigma}_m^2\\
\dotprod{ \bbA}{\bU^1}\\
\vdots \\
\dotprod{ \bbA}{\bU^l}
\end{bmatrix}.
\end{equation}
where $\widehat{\sigma}_m$ is given in Lemma \ref{lemma: estimator of trace of Sigma}.

\end{appendix}

\end{document}